\documentclass{aims}

\usepackage{amsmath}
\usepackage{paralist}
\usepackage[misc]{ifsym}
\usepackage{epsfig}
\usepackage{epstopdf}
\usepackage{amssymb}
\usepackage{algorithm}
\usepackage{algorithmic}
\usepackage[colorlinks=true]{hyperref}
\hypersetup{urlcolor=blue, citecolor=red}
\allowdisplaybreaks
\def\currentvolume{X}
 \def\currentissue{X}
  \def\currentyear{200X}
   \def\currentmonth{XX}
    \def\ppages{X-XX}
     \def\DOI{10.3934/xx.xxxxxxx}

\newcommand{\Ebb}{\mathbb{E}}
\newcommand{\Nbb}{\mathbb{N}}
\newcommand{\Zbb}{\mathbb{Z}}
\newcommand{\Pbb}{\mathbb{P}}
\newcommand{\Rbb}{\mathbb{R}}
\newcommand{\Ubb}{\mathbb{U}}
\newcommand{\Xbb}{\mathbb{X}}
\newcommand{\Ssf}{\mathsf{S}}
\newcommand{\Rsf}{\mathsf{R}}
\newcommand{\Tsf}{\mathsf{T}}

\newcommand{\abs}[1]{\lvert #1\rvert}
\newcommand{\norm}[1]{\lVert #1\rVert}

\DeclareMathOperator{\diag}{diag}
\DeclareMathOperator{\ran}{ran}
\DeclareMathOperator{\rank}{rank}
\DeclareMathOperator{\tr}{tr}

\newtheorem{theorem}{Theorem}[section]
\newtheorem{corollary}[theorem]{Corollary}

\newtheorem{lemma}[theorem]{Lemma}

\theoremstyle{definition}

\newtheorem{remark}[theorem]{Remark}
\newtheorem{assumption}[theorem]{Assumption}

\newtheoremstyle{repeat}
    {\topsep}{\topsep}
    {}{}
    {\bfseries}{.}{ }
    {\thmname{#1}\thmnote{ \bfseries #3}}

\theoremstyle{repeat}
\newtheorem{theorem2}{Theorem}
\newtheorem{corollary2}{Corollary}

\title[Accelerated bistochastic kernel EVD]
{Accelerated decomposition of bistochastic kernel matrices by low rank
approximation}

\author[Chris Vales and Dimitrios Giannakis]{}

\subjclass{Primary: 65F55; Secondary: 65F15.}
\keywords{kernel methods; bistochastic; low rank approximation;
Nystr\"om approximation; partial Cholesky factorization;
spatiotemporal pattern extraction}

\thanks{$^*$Corresponding author: Chris Vales}

\begin{document}
\maketitle

\centerline{\scshape
Chris Vales$^{{\href{mailto:chris.vales@dartmouth.edu}%
    {\textrm{\Letter}}}*1}$
and Dimitrios Giannakis$^{{\href{mailto:dimitrios.giannakis@dartmouth.edu}%
    {\textrm{\Letter}}}1}$}
\medskip

{\footnotesize
 \centerline{$^1$Department of Mathematics, Dartmouth College,
    Hanover, NH, USA}
}
\bigskip

\centerline{(Communicated by Handling Editor)}

\begin{abstract}
We develop an accelerated algorithm for the approximate eigenvalue
decomposition of symmetrically normalized kernel matrices,
focusing on a bistochastic normalization.
Our approach constructs a low rank approximation of the original
kernel matrix by the pivoted partial Cholesky algorithm,
and uses it to compute an approximate decomposition of its
normalization without requiring the formation of the full kernel
matrix.
The cost of the proposed algorithm depends linearly on the size of
the employed training dataset and quadratically on the rank of
the low rank approximation, offering a significant cost reduction
compared to the naive approach.
We derive trace norm error bounds for the approximation of two
classes of normalized kernel matrices.
We apply the proposed algorithm to the kernel based extraction of
spatiotemporal patterns from chaotic Kuramoto--Sivashinsky dynamics,
using a distributed GPU implementation that can scale to large
datasets.
\end{abstract}

\section{Introduction}\label{sec:intro}
Data matrices of large size often arise in the application of
data driven methods to various domains.
An example is given by kernel matrices, whose entries are determined
by the evaluation of a kernel function on a set of data points.
Large kernel matrices arise in kernel based methods such as
support vector machines, applied to tasks such as clustering and
regression \cite{Schoelkopf2001}.
Kernel methods enable the use of linear computational methods
to capture nonlinear relationships in the underlying data.
In addition, they facilitate the application of functional analytic
concepts to datasets without additional mathematical structure,
such as the structure of vector spaces or manifolds.

Despite their large size, data matrices often have a relatively
low approximate rank.
Their underlying low rank structure manifests itself in fast
spectral decay and can be exploited to compute low rank
approximations.
In turn, these approximations can be used to accelerate downstream
computations with only a moderate loss in accuracy.
However, explicitly computing the eigenvalue or singular value
decomposition of a data matrix to construct its low rank
approximation is often unfeasible in practice.
This creates the need for efficient and effective low rank
approximation methods that can scale to matrices of large size
\cite{Halko2011,Martinsson2020,ChenY2025,Epperly2025}.

In this work we are interested in the low rank approximation of
symmetrically normalized kernel matrices, focusing in particular
on a bistochastic normalization.
Using the partial Cholesky factorization method
\cite{ChenY2025},
we build an approximate, low rank factorization of the original
unnormalized kernel matrix.
We then use this factorization to approximate the eigenvalue
decomposition of the considered normalized kernel matrices with
reduced cost, thereby enabling their application to large datasets.
We apply our approach to spatiotemporal pattern extraction
from chaotic Kuramoto--Sivashinsky dynamics using the vector
valued spectral analysis (VSA) method \cite{Giannakis2019vsa}.

The remainder of Section \ref{sec:intro}
is devoted to a review of the symmetric and bistochastic kernel
normalizations considered in this work, and of the pivoted partial
Cholesky algorithm used for their approximation.
In Section \ref{sec:approx} we present the proposed algorithm
for the accelerated eigenvalue decomposition of the normalized
matrices, followed by the derivation of trace norm error bounds
in Section \ref{sec:error}.
In Section \ref{sec:application} we employ the presented algorithm
to extract spatiotemporal patterns from chaotic Kuramoto--Sivashinsky
dynamics in one spatial dimension, followed by a brief conclusion in
Section \ref{sec:conclusion}.
Appendix \ref{app:rpc} provides an overview of the randomly
pivoted Cholesky (RPC) algorithm \cite{ChenY2025,Epperly2025}
employed in this work.
Finally, proofs for the presented error bounds are given in
Appendix \ref{app:proof}.
The code used to generate our numerical results and figures can be
accessed online\footnote{\url{https://github.com/cval26/kernel_evd}}.

\subsection{Normalized kernels}\label{sec:kernel}
We consider the data space $\Ubb=\Rbb^d$, $d\in\Nbb$,
equipped with a probability measure $\mu$
with compact support.
We introduce a kernel function $k$ with the following
properties, which are satisfied by many of the covariance kernels
employed in statistics and machine learning, such as the exponential,
Gaussian and other Mat\'ern kernels
\cite{Schoelkopf2001,Rasmussen2006}.
\begin{assumption}\label{asm:kernel}
    The kernel function $k\colon\Ubb\times\Ubb\to\Rbb$
    is continuous, bounded and positive semidefinite,
    with $k(u,v)>0$ for all $u$, $v\in\Ubb$.
\end{assumption}

Given a collection of $N\in\Nbb$ data samples
$\Ubb_N=\{u_i\}_{i=0}^{N-1}\subset\Ubb$,
we build the kernel matrix
$K\in\Rbb^{N\times N}$,
whose entries $K_{ij}=k(u_i,u_j)$
correspond to evaluation of kernel function $k$
on the available samples.
By construction, the kernel matrix
$K$ is positive semidefinite,
meaning that it is symmetric and that
$y^\top Ky\geq 0$
for every $y\in\Rbb^N$.
As a consequence, it has a well defined eigenvalue decomposition
with real nonnegative eigenvalues and orthogonal eigenvectors.
In what follows, we denote by
$\mu_N=\sum_{i=0}^{N-1}\delta_{u_i}/N$
the empirical sampling measure associated with the
data samples $\Ubb_N$.
In applications we typically only have access to the finite collection
of samples $\Ubb_N$.
As such, we use the sampling measure $\mu_N$ to approximate $\mu$,
assuming an appropriate form of weak convergence as $N\to\infty$.

It is often desirable in applications to normalize the kernel function
$k$ while maintaining its symmetry.
One way to achieve that is to define the normalized kernel
$\ell\colon\Ubb\times\Ubb\to\Rbb$
$$\ell(u,v)=\frac{k(u,v)}{\sqrt{d(u)}\sqrt{d(v)}}$$
with positive normalization function $d\colon\Ubb\to\Rbb$
$$d(u)=\int_\Ubb k(u,v)d\mu(v)$$
whose positivity follows from Assumption \ref{asm:kernel}.
For this choice of normalization, the integral
$\int_\Ubb\ell(u,v)d\mu(v)$
is not necessarily equal to one for every $u\in\Ubb$,
meaning that $\ell$ is generally not a stochastic
(Markov) kernel.
In the discrete setting, $\mu$ is replaced by $\mu_N$
and the above normalization procedure corresponds to forming
the $N\times N$ normalized kernel matrix
$$L=D^{-1/2}KD^{-1/2}$$
with diagonal matrix $D=\diag(K1_N)$
holding the row sums of $K$ in its diagonal,
and $1_N\in\Rbb^N$ denoting the vector of ones.
Normalized kernel matrices of this form are often employed
in diffusion maps algorithms \cite{Coifman2006}
and applications such as kernel spectral clustering
\cite{Ng2001,Zass2005,vonLuxburg2007}.

Normalizing the kernel $k$ in a way that turns it
into a stochastic kernel while maintaining its symmetry can
be achieved with the bistochastic normalization procedure
developed in \cite{Coifman2013}.
The bistochastic kernel function
$p\colon\Ubb\times\Ubb\to\Rbb$
is defined as
\begin{equation}\label{eq:bistoch-kernel}
    p(u,v)=\int_\Ubb\frac{k(u,w)k(w,v)}{d(u)q(w)d(v)}\,d\mu(w)
\end{equation}
with positive functions
$d\colon\Ubb\to\Rbb$ and $q\colon\Ubb\to\Rbb$
$$d(u)=\int_\Ubb k(u,v)d\mu(v)\qquad
    q(u)=\int_\Ubb\frac{k(u,v)}{d(v)}d\mu(v).$$
The analogous discrete procedure involves forming the
$N\times N$ bistochastic kernel matrix
$$P=D^{-1}KQ^{-1}KD^{-1}$$
with diagonal matrices
$$D=\diag(K1_N)\qquad
    Q=\diag(KD^{-1}1_N).$$
Being a stochastic matrix, $P$ can be interpreted as the transition
probability matrix of a Markov chain with $N$ states.
As a result, its eigenvalues are within the interval
$[0,1]\subset\Rbb$,
with its leading eigenvalue being equal to one.
More generally, $p$ induces a self-adjoint Markov operator on
$L^2(\Ubb,\mu)$, which preserves integrals, positive functions and
constant functions \cite{Lasota1994}.

Maintaining symmetry while normalizing a kernel matrix ensures
that the resulting matrix has a real eigenvalue decomposition
with orthogonal eigenvectors.
In the continuous setting, this means that the eigenfunctions of
the associated kernel integral operator define an orthonormal basis
of $L^2(\Ubb,\mu)$.
We are going to use this property to extract spatiotemporal patterns
from a dynamical system in Section \ref{sec:application}.
One motivating reason for working with an ergodic
bistochastic kernel integral operator is that its leading eigenvalue
is equal to one and the leading eigenfunction is constant.
This property is useful for proving convergence of Galerkin
approximation schemes based on the obtained eigenfunctions
\cite{Das2019}.

The bistochastic normalization of kernel matrix $K$
can be achieved with other methods as well, most notably with
the Sinkhorn iterative method \cite{Sinkhorn1964,Sinkhorn1967}
and its symmetrized variant \cite{Zass2005},
which iterate between normalizing the rows and columns of a
matrix until a chosen normalization threshold is met.
The normalized kernel resulting from the Sinkhorn algorithm
cannot be written as an explicit transformation of the original
kernel function $k$.
As a result, to use the Sinkhorn method one has to form the full
kernel matrix $K$, preventing the use of this method in cases where
$K$ does not fit in the available memory.
In contrast, as demonstrated in the following sections, the
fact that the bistochastic kernel \eqref{eq:bistoch-kernel}
is available in closed form enables its manipulation without
the formation of the full kernel matrix.
In addition to their computational differences,
the Sinkhorn normalization and the one defined by
\eqref{eq:bistoch-kernel}
generally lead to different bistochastic kernel functions.
An analytical and computational comparison of the two normalized
kernels is left for future work.

\subsection{Partial Cholesky factorization}\label{sec:cholesky}
The Cholesky factorization of a positive semidefinite matrix
$K\in\Rbb^{N\times N}$
is a decomposition of the form $K=FF^\top$
with lower triangular Cholesky factor $F\in\Rbb^{N\times N}$
\cite{Trefethen1997}.
Because $K$ is only assumed to be positive semidefinite,
instead of positive definite, its Cholesky factorization is
generally not unique.
In this work we focus on the \emph{partial} Cholesky factorization,
which is an approximate decomposition
$$K\approx\tilde K=FF^\top$$
where the partial Cholesky factor
$F\in\Rbb^{N\times r}$
need not be lower triangular anymore.
We will refer to $\tilde K$ as a rank-$r$ partial Cholesky
factorization of $K$, meaning that $F$ has $r$ columns and that
$\rank\tilde K\leq r$ with rank parameter $r<N$.
Although the full and partial Cholesky factorizations can be used
with general positive semidefinite matrices, in this work we focus
specifically on positive semidefinite kernel matrices.

The partial Cholesky factorization can be used to compute low rank
approximations of positive semidefinite matrices by the
pivoted partial Cholesky algorithm
\cite{ChenY2025,Epperly2025}.
Given a rank parameter $r$, the algorithm selects $r$ column
indices (\emph{pivots}) of the input matrix $K$
and uses the corresponding $r$ columns to compute the approximate
factorization $\tilde K=FF^\top$.
The approximation computed by the pivoted partial Cholesky algorithm
for a given set of column indices $\Ssf=\{s_0,\ldots,s_{r-1}\}$
is equal to the column Nystr\"om approximation
$$K\approx\tilde K=K(:,\Ssf)K(\Ssf,\Ssf)^+K(:,\Ssf)^\top$$
where $K(:,\Ssf)$ denotes the $N\times r$ submatrix formed by the
columns indexed by $\Ssf$, $K(\Ssf,\Ssf)$ the $r\times r$
submatrix formed by the rows and columns indexed by $\Ssf$,
and superscript $+$ the Moore--Penrose pseudoinverse
\cite{Drineas2005,ChenY2025}.
Namely, the pivoted partial Cholesky algorithm can be used to
compute the Nystr\"om approximation $\tilde K$ in
the factored form
\begin{equation}\label{eq:nystrom-factored}
\tilde K=FF^\top=\bigl[K(:,\Ssf)L^{-\top}\bigr]
    \bigl[K(:,\Ssf)L^{-\top}\bigr]^\top
\end{equation}
based on the Cholesky decomposition $K(\Ssf,\Ssf)=LL^\top$
with lower triangular Cholesky factor $L\in\Rbb^{r\times r}$.
As we demonstrate in the following sections, having direct
access to either $F$ or $L$ enables the design of efficient
approximation algorithms for normalized kernel matrices.

Throughout this work we employ the following assumption, which
is always satisfied when $\tilde K$ is the column Nystr\"om
approximation of $K$ corresponding to a subset of its columns
\cite{ChenY2025}.
\begin{assumption}\label{asm:psd}
    $K\geq\tilde K\geq 0$
    in the positive semidefinite matrix order.
\end{assumption}

The column indices sampled by the pivoted partial Cholesky algorithm
are usually selected iteratively based on the diagonal entries of the
input matrix $K$, with different sampling strategies
leading to different variations of the algorithm and different
approximation error bounds.
The indices selected by the employed algorithm can be understood
as the most ``important'' states of the dataset according to the
used sampling strategy; they are often referred to as
\emph{landmarks} or \emph{cluster centers}.
Namely, every pivoted partial Cholesky algorithm is in effect
also a sampling algorithm that can be used to subsample a given
dataset in an informed manner.
Importantly, the state sampling and the factorization of the
associated Nystr\"om approximation are carried out simultaneously
by the same algorithm.

In this work we employ the adaptive random sampling strategy
of RPC \cite{ChenY2025},
using the trace norm error $\norm{K-\tilde K}_*=\tr(K-\tilde K)$
to measure the degree of accuracy of the resulting approximation.
For our numerical results in Section \ref{sec:application},
we are going to use the low memory, accelerated version (ARPC-LM)
of the algorithm \cite{Epperly2025},
which samples pivots in blocks instead of one by one,
using an additional accept/reject check to ensure that the
resulting pivot acceptance distribution matches the sampling
distribution of RPC.
In addition, ARPC-LM does not store the $N\times r$ factor matrix $F$,
but works via the smaller $r\times r$ Cholesky factor $L$
introduced in \eqref{eq:nystrom-factored},
thereby enabling its application to large datasets.
Overall, the asymptotic computational cost of computing a rank-$r$
partial Cholesky factorization of an $N\times N$ matrix
$K$ is $O(Nr^2)$.
In Appendix \ref{app:rpc} we provide a brief overview
of ARPC-LM and pseudocode for its implementation.

For an $N\times N$ positive semidefinite matrix $K$
and a rank parameter $r'\in\Nbb$,
we denote by $[K]_{r'}\in\Rbb^{N\times N}$
an optimal rank-$r'$ approximation of $K$,
which is generally not unique.
The approximation $[K]_{r'}=V\Lambda_{r'}V^\top$
can be defined via its truncated eigenvalue decomposition (EVD),
where $\Lambda_{r'}$ denotes the diagonal matrix with the trailing
$N-r'$ eigenvalues set to zero, assuming ordering in descending order.
The approximation is optimal in the sense that, given the rank $r'$,
the error $\norm{K-[K]_{r'}}$ attains its smallest value
in all Schatten matrix norms.
The nonuniqueness may arise in cases where the algebraic multiplicity
of the eigenvalues prevents a unique sorting in descending order.

We state a modified version of \cite[Theorem 2.3]{ChenY2025} below.

\begin{theorem}[{\cite[Theorem 2.3]{ChenY2025}}]\label{thm:rpc}
    Fix $r'\in\Nbb$ and $\epsilon>0$ and let $K$ be a positive
    semidefinite matrix.
    Let $\tilde K$ denote a rank-$r$ partial Cholesky factorization
    of $K$, derived by the randomly pivoted Cholesky algorithm
    using $r\geq r'$ columns of $K$.
    The approximation $\tilde K$ satisfies
    $$\Ebb\tr(K-\tilde K)\leq\zeta,\qquad
        \zeta=(1+\epsilon)\tr(K-[K]_{r'})$$
    for every $r$ such that
    $$r\geq\frac{r'}{\epsilon}+r'\log\frac{1}{\epsilon\eta},\qquad
        \eta=\tr(K-[K]_{r'})/\tr K$$
    where the expectation is taken over the random choice
    of sampled columns.
\end{theorem}

\begin{remark}\label{rem:rpc}
\ 
\begin{enumerate}
\item The theorem shows that when the number of sampled columns
    $r$ of the partial Cholesky factorization satisfies the presented
    inequality, then the expected trace norm error of the approximation
    is at most a factor $1+\epsilon$ larger than the best error attainable
    for an approximation of lower rank $r'\leq r$.
    As a result, the error bound $\zeta$ can also be written in the form
    $$\zeta=(1+\epsilon)\sum_{i=r'}^{N-1}\lambda_i$$
    using the trailing $N-r'$ eigenvalues $\lambda_i$ of $K$
    that are not captured by $[K]_{r'}$.
\item In the accelerated algorithm ARPC-LM used for our numerical
    results, the pivots are sampled in blocks instead of one by one.
    Epperly et al \cite{Epperly2025} have extended
    Theorem \ref{thm:rpc} to use the number of blocks of pivots
    required to ensure a certain level of accuracy, instead of the
    number of individual pivots $r$.
    For the error analysis presented in Section \ref{sec:error},
    we use the version of the theorem as stated above.
    The results are independent of whether the number of blocks
    or the number of individual pivots is used to control the
    expected error $\zeta$.
\end{enumerate}
\end{remark}

\section{Normalized kernel approximation}\label{sec:approx}
There is a large body of literature on the use of
subsampling and low rank approximation algorithms to accelerate
the implementation of kernel methods and enable their application
to large datasets
\cite{Williams2000,Fowlkes2001,deSilva2004,Fowlkes2004,Drineas2005,
Belabbas2009pnas,Belabbas2009rsa,Kumar2009,LiM2011,Kumar2012,
Choromanska2013,Gittens2013,Lau2020,Pourkamali2020,ChenY2025}.
In this work, we focus on applications that require the
computation of the EVD of a normalized kernel matrix,
with kernel spectral clustering being a prime example
\cite{Ng2001,Zelnik2004,Zass2005,vonLuxburg2007}.

We begin with a unifying overview of the approximation method
considered in this work, before specializing it to the symmetric
and bistochastic normalized kernel matrices.
Motivated by the partial Cholesky factorization reviewed in
Section \ref{sec:cholesky},
we consider the approximation of the EVD of a
normalized kernel matrix $M\in\Rbb^{N\times N}$
that can be written in the approximate factorized form
$$M\approx\tilde M=BCB^\top$$
with $B\in\Rbb^{N\times r}$ and symmetric $C\in\Rbb^{r\times r}$,
where $r<N$ is a chosen rank parameter.
The important aspect of this factorization is that it is rank
revealing, since $B$ is of reduced dimension $N\times r$.
Using the above factorization, the EVD of $\tilde M$
can be computed in two main steps.
First, we compute the reduced QR factorization $B=YR$
to rewrite
$$\tilde M=Y(RCR^\top)Y^\top$$
with $Y\in\Rbb^{N\times r}$ and $R\in\Rbb^{r\times r}$.
Second, we compute the EVD of the small $r\times r$ matrix
$RCR^\top=V\Lambda V^\top$,
with eigenvector matrix $V\in\Rbb^{r\times r}$
and diagonal eigenvalue matrix $\Lambda\in\Rbb^{r\times r}$.
This leads to the EVD
$$\tilde M=U\Lambda U^\top$$
with the orthonormal columns of $U=YV\in\Rbb^{N\times r}$
holding the eigenvectors.
This general method involves an asymptotic computational cost
$O(Nr^2)$, dominated by the cost of the QR factorization.
The method computes the exact EVD of the approximate
matrix $\tilde M$, whose accuracy will be studied analytically
in Section \ref{sec:error}.

In Section \ref{sec:approx-symm} we begin by demonstrating this
method with the symmetric normalization $L$.
Next, in Section \ref{sec:approx-bistoch} we extend our approach to
the bistochastic normalization $P$.
In Sections \ref{sec:nystrom} and \ref{sec:alternatives}
we discuss connections with the Nystr\"om method and
alternative approximation strategies.

\subsection{Symmetric normalization}\label{sec:approx-symm}
We consider a kernel function $k$ satisfying
Assumption \ref{asm:kernel}.
Using a collection of $N\in\Nbb$ state samples
$\Ubb_N=\{u_i\}_{i=0}^{N-1}\subset\Ubb$,
we build the corresponding $N\times N$ kernel matrices
$K_{ij}=k(u_i,u_j)$
and
$$L=D^{-1/2}KD^{-1/2}$$
with diagonal matrix $D=\diag(K1_N)$
with positive diagonal entries.
Our goal is to compute the EVD of matrix $L$,
which is an operation of asymptotic cost $O(N^3)$.
We operate under the assumption that computing the EVD of $L$
directly is not feasible due to the large number of samples $N$.
For this reason, we are going to perform a low rank approximation
of $L$ and use it to reduce the cost of the eigenvalue problem.

We employ a pivoted partial Cholesky algorithm to build an
approximation of matrix $K$ that is of lower rank $r<N$
$$K\approx\tilde K=FF^\top$$
with partial Cholesky factor $F\in\Rbb^{N\times r}$.
Next, we build the $N\times N$ diagonal matrix
$\tilde D=\diag(\tilde K1_N)$
and form the normalized low rank approximation
$$L\approx\tilde L=\tilde D^{-1/2}\tilde K\tilde D^{-1/2}
    =(\tilde D^{-1/2}F)(\tilde D^{-1/2}F)^\top.$$
To form $\tilde L$ we assume that the diagonal entries of
$\tilde D$ are positive, leading to a well defined normalization.
This can be ensured by picking a sufficiently large rank parameter $r$:
the difference $K-\tilde K$ approaches zero in the trace norm as
$r\to N$, which implies that it approaches zero entrywise.
In turn, this implies that the row sums of $\tilde K$ approach
those of $K$, which are positive by Assumption \ref{asm:kernel}.
Estimates of this form are given
in Lemma \ref{lem:dbound}
in Appendix \ref{app:proof}.

The approximate matrix $\tilde L$ is now in the rank revealing
form $BCB^\top$ introduced earlier,
with $N\times r$ matrix $B=\tilde D^{-1/2}F$
and the $r\times r$ identity matrix $C=I$.
Computing the EVD of $\tilde L$ requires the reduced QR
factorization $\tilde D^{-1/2}F=YR$
and the EVD $RR^\top=V\Lambda V^\top$,
leading to the desired EVD
$$\tilde L=(YV)\Lambda(YV)^\top.$$

For this class of normalized kernel matrices, the inner matrix
$C$ is equal to the identity matrix.
As a result, the EVD of $\tilde L$ can also be obtained by
computing the reduced singular value decomposition (SVD) of
$\tilde D^{-1/2}F$.
Both method variations involve a cost $O(Nr^2)$ and lead to
the same decomposition up to sign differences in the eigenvectors.
The accuracy of the SVD-based variant in kernel spectral clustering
has been demonstrated in \cite{ChenY2025},
using RPC to build the partial Cholesky factor $F$.

\subsection{Bistochastic normalization}\label{sec:approx-bistoch}
We now focus on the bistochastic normalized kernels outlined
in Section \ref{sec:kernel}.
As before, we use a kernel function $k$ satisfying
Assumption \ref{asm:kernel}
and a collection of $N$ state samples $\Ubb_N\subset\Ubb$.
We build the $N\times N$ kernel matrix $K_{ij}=k(u_i,u_j)$
and its bistochastic normalization
$$P=D^{-1}KQ^{-1}KD^{-1}$$
with $N\times N$ diagonal matrices
$D=\diag(K1_N)$ and $Q=\diag(KD^{-1}1_N)$
with positive diagonal entries.
As earlier, we assume that computing the EVD of $P$ directly
is not feasible due to the large number of samples $N$.

Using a pivoted partial Cholesky algorithm, we compute
a rank-$r$ approximation $\tilde K$
$$K\approx\tilde K=FF^\top$$
with partial Cholesky factor
$F\in\Rbb^{N\times r}$, $r< N$.
Assuming again that $r$ is sufficiently large so that the
normalization matrices are well defined
$$\tilde D=\diag(\tilde K1_N)\qquad
    \tilde Q=\diag(\tilde K\tilde D^{-1}1_N)$$
we form the low rank approximation of $P$
$$\tilde P=\tilde D^{-1}\tilde K\tilde Q^{-1}\tilde K\tilde D^{-1}
    =(\tilde D^{-1}F)(F^\top\tilde Q^{-1}F)(\tilde D^{-1}F)^\top.$$

The approximate matrix $\tilde P$ is now expressed in the rank
revealing form $BCB^\top$ introduced earlier,
with $N\times r$ matrix $B=\tilde D^{-1}F$
and $r\times r$ symmetric matrix $C=F^\top\tilde Q^{-1}F$.
As such, computing the EVD of $\tilde P$ requires the reduced
QR factorization $\tilde D^{-1}F=YR$
and the EVD
$RF^\top\tilde Q^{-1}FR^\top=V\Lambda V^\top$,
leading to the desired EVD
$$\tilde P=(YV)\Lambda(YV)^\top$$
with eigenvector matrix $U=YV$.

Pseudocode for the implementation of the method is given in
Algorithm \ref{alg:bistoch}.
The algorithm takes as inputs the kernel matrix $K$,
approximation rank $r$ and block size $b$ needed for ARPC-LM,
and returns the $r$ leading eigenvectors $U$
and eigenvalues $\Lambda$ of the approximation $\tilde P$.
To be more precise, the algorithm returns at most $r$ leading
eigenvectors and associated eigenvalues, due to the way the
ARPC-LM algorithm constructs the matrix $L^{-1}$; see
Appendix \ref{app:rpc} for details.
Importantly, the algorithm never forms the full kernel matrices
$K$ or $P$.
As presented, Algorithm \ref{alg:bistoch}
uses the CholeskyQR method to compute the QR decomposition
of the tall $N\times r$ matrix $\tilde D^{-1}F$,
and includes a regularization term $\gamma>0$
for the Cholesky step \cite{Yamamoto2015}.
The algorithm involves an asymptotic computational cost
$O(Nr^2)$ for the rank-$r$ partial Cholesky factorization,
$O(r^3)$ for the EVD of the inner $r\times r$ matrix,
and $O(Nr^2)$ for the QR factorization and
some of the intermediate matrix products,
yielding an overall asymptotic cost $O(Nr^2)$.

\begin{algorithm}[t]
    \caption{Approximate EVD for bistochastic kernels}
\label{alg:bistoch}
\begin{algorithmic}
    \STATE Input: $N\times N$ kernel matrix $K$,
        approximation rank $r<N$, block size $b$
    \STATE Output: $N\times r$ eigenvector matrix $U$,
        $r\times r$ diagonal eigenvalue matrix $\Lambda$
    \STATE
    \STATE $\Ssf,\,L^{-1}\gets\texttt{arpclm}(K, r, b)$
        \hfill ARPC-LM (Alg.\@ \ref{alg:arpclm})
    \STATE $F\gets K(:,\Ssf)L^{-\top}$\hfill
        Eq.\@ \eqref{eq:nystrom-factored}
    \STATE
    \STATE $\tilde D\gets\texttt{diag}(FF^\top 1_N)$
        \hfill diagonal normalization matrix
    \STATE \texttt{assert} $\texttt{diag}(\tilde D)>0$
        \hfill ensure positivity
    \STATE $\tilde Q\gets\texttt{diag}(FF^\top\tilde D^{-1}1_N)$
        \hfill diagonal normalization matrix
    \STATE \texttt{assert} $\texttt{diag}(\tilde Q)>0$
        \hfill ensure positivity
    \STATE
    \STATE $R\gets\texttt{chol}\bigl((\tilde D^{-1}F)^\top\tilde D^{-1}F
        +\gamma I\bigr)$\hfill $r\times r$ regularized Cholesky
    \STATE $V,\,\Lambda\gets
        \texttt{evd}(RF^\top\tilde Q^{-1}FR^\top)$
        \hfill $r\times r$ EVD, descending order
    \STATE $Y\gets \tilde D^{-1}FR^{-1}$\hfill
        such that $\tilde D^{-1}F=YR$
    \STATE $U\gets YV$\hfill $N\times r$ eigenvector matrix
\end{algorithmic}
\end{algorithm}

Our results in Section \ref{sec:application}
are generated using a distributed GPU implementation of
Algorithm \ref{alg:bistoch},
similar in design to the implementation of
Algorithm \ref{alg:arpclm} described in
Appendix \ref{app:rpc}.
The implementation of Algorithm \ref{alg:bistoch}
distributes all computations and the storage of the
required $r\times r$ matrices among the available GPUs.
In addition, $N\times r$ matrices such as
$K(:,\Ssf)$, $F$ and $Y$ are never formed in memory;
rather, their values are regenerated as needed in blocks of
rows, with the size of the regeneration block set by the
amount of memory available.
Excluding $N$-sized vectors, this renders the memory cost of
Algorithm \ref{alg:bistoch}
independent of $N$, enabling its application to large datasets.
As with the implementation of Algorithm \ref{alg:arpclm},
all $N$-dependent computations are performed in single precision
(32-bit) floating point arithmetic, whereas the computations
involving $r\times r$ matrices are performed in double precision
(64-bit).
Finally, the implemented algorithm allows the extraction of only
a subset of the leading $r$ eigenvectors (columns of $U$),
removing the need to form the full matrix $U$.
Additional details on the implementation can be found in the accompanying
code\footnote{\url{https://github.com/cval26/kernel_evd}}.

\subsection{Nystr\"om method connection}\label{sec:nystrom}
As explained in Section \ref{sec:cholesky}, the partial Cholesky
factorization $\tilde K=FF^\top$ corresponds to
a factored form of the column Nystr\"om approximation
$\tilde K=K(:,\Ssf)K(\Ssf,\Ssf)^+K(\Ssf,:)$
associated with the set of sampled indices $\Ssf$.
As a result, the EVD approximation method presented
above is closely connected to the Nystr\"om approximation method,
where one uses the matrices $K(:,\Ssf)$ and $K(\Ssf,\Ssf)$,
instead of the factor matrix $F$.

In what follows, we review the connection to the Nystr\"om method
using the unnormalized kernel matrix $K$, but a similar analysis
can be performed for normalized kernel matrices
\cite{Choromanska2013,Pourkamali2020}.
Assuming access to a set of sampled indices $\Ssf$ of size $r$,
we use it to construct the $r\times r$ kernel matrix
$K(\Ssf,\Ssf)$ and compute its EVD directly,
$K(\Ssf,\Ssf)=VZV^\top$.
The computed eigenvectors $V$ can be extended to approximate
the eigenvectors of the original matrix $K$
by using the cross-kernel matrix $K(:,\Ssf)$
$$\tilde K=(K(:,\Ssf)V)Z^+(K(:,\Ssf)V)^\top
    =\tilde VZ^+\tilde V^\top.$$
The extension operation $K(:,\Ssf)VZ^+$ is the discrete analog
of the Nystr\"om extension formula for kernel integral operators
written in matrix form \cite{Drineas2005}.

The extended eigenvectors $\tilde V$ can be used as an approximation
of the eigenvectors of the original kernel matrix $K$.
However, they do not correspond to the eigenvectors of the
approximation $\tilde K$; since the columns of $\tilde V$ are
generally not orthogonal, the above decomposition of $\tilde K$
in terms of $\tilde V$ is generally not an EVD.
One way of turning the above into the EVD of $\tilde K$ is
the following sequence of steps: compute the QR factorization
$\tilde V=YR$, followed by the EVD of the inner matrix
$RZ^+R^\top=W\Lambda W^\top$.
This results in the EVD
$$\tilde K=(YW)\Lambda (YW)^\top=U\Lambda U^\top$$
and mirrors precisely the two main steps involved in the
general EVD approximation algorithm presented in
Section \ref{sec:approx}.
Given the uniqueness of the EVD of $\tilde K$,
the two approaches arrive at the same result, up to changes in
sign in the eigenvectors.
Conceptually, their main difference is that one uses the factor
matrix $F$, whereas the other works with the matrices
$K(\Ssf,\Ssf)$ and $K(:,\Ssf)$.

Although theoretically equivalent, we argue that the approach
that relies on the factor matrix $F$ is better suited to the
computation of the EVD of $\tilde K$.
This is because it avoids the need to first decompose
the small kernel matrix $K(\Ssf,\Ssf)$, extend this decomposition
to $\tilde K$, and finally turn that into the EVD of $\tilde K$.
Rather, using $F$ allows us to move directly to the EVD of
$\tilde K$ in two main steps.
This becomes especially useful when normalized kernel matrices
are involved, such as the matrices $L$ and $P$ studied in this work.
Finally, using the partial Cholesky factorization to represent
the Nystr\"om approximation in factored form enables the combination
of the sampling and kernel evaluation steps into one algorithm.
This is the most natural approach for adaptive sampling strategies,
such as the sampling strategy of RPC \cite{ChenY2025,Epperly2025}
used in the present work.

\subsection{Alternative approximation strategies}
\label{sec:alternatives}
As mentioned earlier, computing the EVD of kernel matrix $K$
directly carries an asymptotic cost $O(N^3)$ and requires
the formation of the full matrix in memory, which is often
unfeasible in practice.
Although a partial Cholesky factorization of $K$ is employed to
alleviate this cost in this work, there are several alternative
methods that can also be considered.

Iterative Krylov methods can be used to approximate the
leading eigenvalues and associated eigenvectors of $K$, without
performing any approximation of the matrix
\cite{Trefethen1997}.
Although they do not require forming $K$ in memory, these methods
involve a significant number of matrix-vector products, meaning that
the full matrix $K$ will have to be recomputed several times.
Randomized variants of these methods are also available, coupling
the iterative Krylov approach with a randomized approximation
of $K$ \cite{Tropp2023}.

Alternative low rank approximation methods include the use of the
randomized SVD method, specifically its Nystr\"om variant
tailored to positive semidefinite matrices
\cite[Algorithms 5.5 \& 5.6]{Halko2011}
\cite[Algorithm 5.6]{Tropp2023}.
These methods employ a matrix of random entries to build an approximate
orthonormal basis for the range space of $K$, then use that to compute
the column Nystr\"om approximation associated with the constructed
basis.
Their implementation requires the multiplication of the kernel matrix
$K$ with a random matrix to build the required basis.
Their asymptotic cost is $O(N^2r)$ and is dominated by the
matrix-matrix multiplication of $K$ with the random matrix used
to approximate its range.
Excluding this multiplication, the cost of the rest of the algorithm
is $O(Nr^2)$, matching the partial Cholesky methods.
Compared to the randomized Nystr\"om methods that use a random
matrix to approximate the range of $K$, partial Cholesky
methods---randomized or not---take advantage of the fact that
$K$ is a kernel matrix and try to identify cluster centers in the
underlying dataset that can be used to approximate $K$.
This allows them to avoid computing all entries of $K$ or
its product with another matrix, leading to their favorable
cost scaling.
As such, the partial Cholesky methods can be viewed as a
specialization of the randomized Nystr\"om methods to kernel
matrices.

Another approach is to approximate the kernel matrix $K$ by a
sparse matrix using a predefined sparsity pattern or nearest
neighbors search.
Contrary to the spectral approximation methods outlined above,
this approach does not try to exploit eigenvalue decay but
sparsity or structure in the matrix entries
\cite{Le2013,Wilson2015,Schaefer2021}.
One example of this approach is the sparse Cholesky factorization
\cite{Schaefer2021,ChenY2025sparse},
which involves an asymptotic cost $O(N\log^d(N/\epsilon))$
for an approximation error $\epsilon$ in an appropriate norm. 
The low rank and sparse or structured approximation strategies
are not mutually exclusive; they can be combined to further
accelerate or improve the accuracy of approximation algorithms
\cite{Zhao2024,Kaminetz2026}.
For instance, by exploiting sparsity in the considered matrix $K$,
the randomized Nystr\"om methods mentioned above can be implemented
with cost $O(Nnr+Nr^2)$, where $n$ denotes the average number of
nonzero entries in each row of $K$.
Which method to choose is usually informed by the type of
structure that can be exploited in the problem at hand.

\section{Error bounds}\label{sec:error}
In this section we present bounds for the trace norm error of
the approximation $L-\tilde L$ for the symmetric normalization
and $P-\tilde P$ for the bistochastic one.

In each case, we start by deriving an error bound that holds for
every realization of the generally random matrix $\tilde K$,
where the source of randomness is the random sampling strategy
that may be used to sample columns of the full kernel matrix $K$.
Then, we specialize this result to a high probability bound for
the case where $\tilde K$ is formed using the RPC algorithm,
employing the expectation bound of Theorem \ref{thm:rpc}.

The EVD approximation algorithm presented in
Section \ref{sec:approx}
computes the exact EVD of the approximate kernel matrices
$\tilde L$ and $\tilde P$.
As a result, the trace norm error bounds presented below
are directly applicable to the produced EVD.
In particular, the results place an upper bound on the error in
terms of a number of the trailing eigenvalues $\lambda_i$ of $K$.
This is because Theorem \ref{thm:rpc} connects the trailing eigenvalues
of $K$ to the error of the approximation of $K$ by $\tilde K$,
which is then employed in forming the normalized matrices
$\tilde L$ and $\tilde P$.
For example, the error bound $\norm{P-\tilde P}_*\leq c\,\zeta$
for a positive constant $c$ can be written in the form
$$\norm{P-\tilde P}_*\leq c\,(1+\epsilon)\sum_{i=r'}^{N-1}\lambda_i$$
in the notation of Theorem \ref{thm:rpc} and Remark \ref{rem:rpc}.

We note that Assumptions \ref{asm:kernel} and \ref{asm:psd}
are in effect throughout this section.
The proofs of all results that follow are given in
Appendix \ref{app:proof}.
We use $\norm{\cdot}$ to denote the operator (spectral) norm.

\subsection{Symmetric normalization}\label{sec:error-symm}
We consider the $N\times N$ matrices
$$L=D^{-1/2}KD^{-1/2}\qquad
    \tilde L=\tilde D^{-1/2}\tilde K\tilde D^{-1/2}$$
from Section \ref{sec:approx-symm}.
In addition, we define
$$d_{\min}=\min_{0\leq i<N}d_i\qquad
    \tilde d_{\min}=\min_{0\leq i<N}\tilde d_i\qquad
    \delta=\min(d_{\min},\tilde d_{\min})$$
where $d_i$ and $\tilde d_i$ denote the respective diagonal entries
of $D$ and $\tilde D$.
From Assumption \ref{asm:kernel} it follows that $d_{\min}>0$.
To ensure the positivity of $\tilde d_i$ we make the
following additional assumption.
The fact that the assumption can be satisfied by choosing a
sufficiently large rank parameter $r$ follows from the estimate
in Lemma \ref{lem:dbound} and the positivity of $d_i$.
\begin{assumption}\label{asm:dpos}
    The rank parameter $r$ is large enough such that
    $\tilde d_i=(\tilde K1_N)_i>0$
    for all $i\in\{0,\ldots,N-1\}$.
\end{assumption}

We start by proving the following result, which bounds the trace norm
error of $L-\tilde L$ by that of the approximation $K-\tilde K$,
and holds for any realization of the generally random matrix
$\tilde K$.
\begin{theorem}\label{thm:symm-error1}
    Under Assumption \ref{asm:dpos}, the normalized matrices
    $L$ and $\tilde L$ satisfy
    $$\lVert L-\tilde L\rVert_*\leq\Bigl[\frac{1}{d_{\min}}
        +\frac{\sqrt{N}}{\delta^2}\tr(K)
        +\frac{N}{4\delta^3}\tr(K)\tr(K-\tilde K)\Bigr]
        \tr(K-\tilde K).$$
\end{theorem}

Pushing the error $\tr(K-\tilde K)$ forward under the
symmetric normalization yields a dependence of $\tr(L-\tilde L)$
on the normalization factors (row sums) $d_i$ and $\tilde d_i$.
The factors appearing in the bound of
Theorem \ref{thm:symm-error1}
are controlled by the smallest values
$d_{\min}$ and $\tilde d_{\min}$.

\begin{remark}\label{rem:symm-error1}
\ 
\begin{enumerate}
\item If the normalized matrix $\tilde L$ is formed
    using $D$ instead of $\tilde D$, then only the first
    term of the bound in Theorem \ref{thm:symm-error1}
    will be present
    $$\lVert L-\tilde L\rVert_*\leq\frac{\tr(K-\tilde K)}{d_{\min}}.$$
    Namely, the error propagation is then controlled by
    $d_{\min}$ alone, since both matrices $L$ and $\tilde L$ use $D$
    for their normalization.
    In addition, in this case the matrix $L-\tilde L$ is positive
    semidefinite, which means that
    $\tr(L-\tilde L)=\norm{L-\tilde L}_*$,
    in contrast to the general case where
    $\tr(L-\tilde L)\leq\norm{L-\tilde L}_*$.
    This normalization choice is often referred to as
    \emph{degree-preserving normalization}.
    Of course, using this normalization requires that we form $D$,
    which in turn requires accessing the full kernel matrix $K$.
\item The terms depending on $N$ appear because we have to bound the
    diagonal entries of matrix $D^{-1/2}-\tilde D^{-1/2}$,
    which is in turn used to bound $L-\tilde L$.
    The fact that the row sums $d_i$ and $\tilde d_i$ depend on the
    entries of a full row of $K$ and $\tilde K$ then leads to the
    dependence on $N$ (Appendix \ref{app:proof}).
    These are the terms that become zero when the degree-preserving
    normalization is used.
\item Pourkamali-Anaraki \cite[Theorem 2]{Pourkamali2020}
    derives a similar error bound for $L-\tilde L$
    in the operator norm instead of the trace norm.
    Their proof employs the Taylor expansion of a matrix function,
    and requires an additional smallness hypothesis dictated by the
    radius of convergence of the Taylor expansion.
\end{enumerate}
\end{remark}

Next, we specialize Theorem \ref{thm:symm-error1} to the case where
the approximation $\tilde K$ is built using RPC.
This means that the error $\tr(K-\tilde K)$ can now be bounded in
expectation by invoking Theorem \ref{thm:rpc}.
Given that $\tilde K$ is a random variable,
$\tilde d_{\min}$ and $\delta$ are also random and correlated
with $\tilde K$.
As a result, to use the error bound of Theorem \ref{thm:rpc}
we condition on a high probability event and prove
the following result.
\begin{corollary}\label{cor:symm-error2}
    Under the assumptions of Theorem \ref{thm:rpc},
    choose $t\in\Rbb$ and the rank parameter $r$ such that
    \begin{equation}\label{eq:t-admit}
        1<t<\frac{d_{\min}}{\sqrt{N}\zeta}.
    \end{equation}
    Then with probability at least $1-1/t$
    $$\lVert L-\tilde L\rVert_*<\Bigl[\frac{1}{d_{\min}}
        +\frac{\sqrt{N}\tr(K)}{(d_{\min}-\sqrt{N}t\zeta)^2}
        +\frac{N\tr(K)t\zeta}{4(d_{\min}-\sqrt{N}t\zeta)^3}\Bigr]t\zeta$$
    with $\tilde d_{\min}>d_{\min}-\sqrt{N}t\zeta>0$.
\end{corollary}

\begin{remark}\label{rem:symm-error2}
\ 
\begin{enumerate}
\item On the high probability event used in
    Corollary \ref{cor:symm-error2}, the positivity of the row sums
    of $\tilde K$ is ensured by the positivity of the row sums of $K$
    and the condition \eqref{eq:t-admit}.
\item In the degree-preserving normalization case, one can
    take the expected value of the bound of Theorem
    \ref{thm:symm-error1} directly, since $d_{\min}$
     is not a random variable.
    Nevertheless, the high probability bound of Corollary
    \ref{cor:symm-error2} is also valid, with only the first
    factor ($1/d_{\min}$) present.
\item The feasibility of condition \eqref{eq:t-admit}
    can be ensured by choosing a sufficiently large rank
    parameter $r$.
    More specifically, the condition is feasible when
    $\sqrt{N}\zeta<d_{\min}$, where $\zeta$ denotes the expected
    error bound from Theorem \ref{thm:rpc}.
    Namely, given the size of the dataset $N$ and the minimum
    row sum $d_{\min}$, one has to choose $\zeta$ small enough to
    ensure a reasonable window of admissible values for $t$.
    The error bound $\zeta$ is in turn controlled by the number
    of columns $r$ of $K$ sampled to produce the approximation
    $\tilde K$, with larger $r$ generally leading to smaller $\zeta$.
    The $\sqrt{N}$ factor enters because we use the trace bound of
    $K-\tilde K$ to place a bound on the difference of row sums
    $d_i-\tilde d_i$ (Appendix \ref{app:proof}).
\item Choosing the value of $t$ requires trading the failure
    probability $1/t$ against the tightness of the trace norm
    error bound.
    As $t\to 1^+$, $d_{\min}-\sqrt{N}t\zeta$ approaches its largest
    admissible value, which tightens the bound of $\norm{L-\tilde L}_*$,
    but the failure probability approaches $1$.
    In the other extreme, as $t\to d_{\min}/(\sqrt{N}\zeta)$
    from below, the failure probability approaches its smallest value,
    but the denominator $d_{\min}-\sqrt{N}t\zeta$ approaches $0$,
    which makes the error bound uninformative.
    Generally, we want the admissible window of $t$ to be large enough
    to allow us to choose a value $t$ that strikes a good compromise
    between a small failure probability and a tight error bound.
\item In addition to $r$, the expected error bound $\zeta$ can also be
    controlled by the kernel function $k$ used to produce the kernel
    matrix $K$.
    For example, for the case of a Gaussian kernel $k$ with
    bandwidth $\varepsilon>0$, increasing the bandwidth generally
    reduces the approximate rank of the resulting $K$.
    This is because the kernel function averages the underlying dataset
    more aggressively, producing a $K$ with fast spectral decay but
    also failing to capture fine features of the data.
    In turn, this means that a smaller error bound $\zeta$
    can be achieved with a relatively modest rank parameter $r$.
    In addition, increasing the bandwidth $\varepsilon$ generally
    increases $d_{\min}$, which also leads to a larger window of
    admissible values for $t$.
\end{enumerate}
\end{remark}

\subsection{Bistochastic normalization}\label{sec:error-bistoch}
We consider the $N\times N$ matrices
$$P=D^{-1}KQ^{-1}KD^{-1}\qquad
    \tilde P=\tilde D^{-1}\tilde K\tilde Q^{-1}\tilde K\tilde D^{-1}$$
from Section \ref{sec:approx-bistoch}
and define the quantities
$$d_{\min}=\min_{0\leq i<N}d_i\qquad
    \tilde d_{\min}=\min_{0\leq i<N}\tilde d_i\qquad
    \delta_d=\min(d_{\min},\tilde d_{\min})$$
and
$$q_{\min}=\min_{0\leq i<N}q_i\qquad
    \tilde q_{\min}=\min_{0\leq i<N}\tilde q_i\qquad
    \delta_q=\min(q_{\min},\tilde q_{\min})$$
with $d_i$ and $\tilde d_i$ the respective diagonal entries of
$D$ and $\tilde D$, and $q_i$ and $\tilde q_i$ those of
$Q$ and $\tilde Q$.
It follows from Assumption \ref{asm:kernel} that
$d_{\min}>0$ and $q_{\min}>0$,
and from Assumption \ref{asm:dpos} that $\tilde d_{\min}>0$.
To ensure the positivity of $\tilde q_{\min}$ we make the
following additional assumption, which can be ensured by
choosing a sufficiently large rank parameter $r$.
This follows from the positivity of $q_i$ and the estimate
in Lemma \ref{lem:qbound}.
\begin{assumption}\label{asm:qpos}
    The rank parameter $r$ is large enough such that
    $\tilde q_i=(\tilde K\tilde D^{-1}1_N)_i>0$
    for all $i\in\{0,\ldots,N-1\}$.
\end{assumption}

The introduced matrices can be written in the Gram form
$P=G^\top G$ and $\tilde P=\tilde G^\top\tilde G$
with
\begin{equation}\label{eq:gram}
G=Q^{-1/2}KD^{-1}\qquad
    \tilde G=\tilde Q^{-1/2}\tilde K\tilde D^{-1}.
\end{equation}
We use the Gram form of $P$ and $\tilde P$ to prove the following
error bound, which holds for any realization of the generally random
matrix $\tilde K$.
\begin{theorem}\label{thm:bistoch-error1}
    Under Assumptions \ref{asm:dpos} and \ref{asm:qpos},
    the normalized matrices $P$ and $\tilde P$ satisfy
    $$\norm{P-\tilde P}_*\leq(1+\norm{\tilde G})\Bigl[
        \frac{1}{\sqrt{q_{\min}}d_{\min}}
        +\frac{\sqrt{N}\tr(K)}{\sqrt{\delta_q}\delta_d^2}
        +\frac{\sqrt{N}\tr(K)}{2\delta_q^{3/2}\delta_d^2}
        \bigl(1+\frac{\norm{K}}{\delta_d}\bigr)\Bigr]
        \tr(K-\tilde K).$$
\end{theorem}
In analogy to Theorem \ref{thm:symm-error1} presented earlier,
the error bound of Theorem \ref{thm:bistoch-error1}
depends on the smallest normalization constants
$d_{\min}$, $\tilde d_{\min}$
and $q_{\min}$, $\tilde q_{\min}$.

\begin{remark}\label{rem:bistoch-error1}
\ 
\begin{enumerate}
\item In analogy to the degree-preserving normalization of
    Remark \ref{rem:symm-error1}, the normalized matrix $\tilde P$
    can be formed using matrices $D$ and $Q$ instead of the
    approximate $\tilde D$ and $\tilde Q$; we refer to this
    as the \emph{structure-preserving normalization}.
    In this case, the error bound simplifies to
    $$\norm{P-\tilde P}_*\leq(1+\norm{\tilde G})
        \Bigl[\frac{1}{\sqrt{q_{\min}}d_{\min}}\Bigr]
        \tr(K-\tilde K)$$
    with the error propagation controlled by $d_{\min}$
    and $q_{\min}$ only.
    In addition to requiring access to the full kernel matrix
    $K$, the structure-preserving normalization does not ensure
    that the rows and columns of $\tilde P$ sum to 1.
\item If one assumes that the entries of $\tilde K$ are nonnegative,
    then the same is true for $\tilde P$ and $\tilde G$.
    This implies that $\norm{\tilde G}=1$, leading to the
    simplified prefactor $1+\norm{\tilde G}=2$.
\item Similarly to Remark \ref{rem:symm-error1}, the bound terms
    depending on $N$ are present because we have to bound the
    diagonal entries of $D^{-1}-\tilde D^{-1}$
    and $Q^{-1/2}-\tilde Q^{-1/2}$ (Appendix \ref{app:proof}).
\item To prove Theorem \ref{thm:bistoch-error1}
    we place the trace norm on the matrices $K$ and $\tilde K$,
    and the operator norm on terms including the normalization matrices
    $D^{-1}$, $\tilde D^{-1}$ and $Q^{-1/2}$, $\tilde Q^{-1/2}$.
    An alternative is to do the opposite, placing the trace norm on the
    normalization matrices instead, in which case the trace norms of
    $D^{-1}-\tilde D^{-1}$ and $Q^{-1/2}-\tilde Q^{-1/2}$
    are bounded by the sums of $\abs{d_i-\tilde d_i}$
    and $\abs{q_i-\tilde q_i}$ instead of their maximum values.
    The choice made in Theorem \ref{thm:bistoch-error1}
    is generally tighter when $\norm{K}_*$ is close to $\norm{K}$,
    which corresponds to fast eigenvalue decay for $K$ and is
    the case where a low rank approximation can be especially
    effective.
\end{enumerate}
\end{remark}

We now specialize Theorem \ref{thm:bistoch-error1}
to the case where the approximation $\tilde K$ is built using RPC,
which enables the use of the expectation bound of
Theorem \ref{thm:rpc} for $\tr(K-\tilde K)$.
As in the previous section, we do that by conditioning on an
appropriate, high probability event.
\begin{corollary}\label{cor:bistoch-error2}
    Under the assumptions of Theorem \ref{thm:rpc},
    let $t\in\Rbb$ and
    $$\alpha_t=d_{\min}-\sqrt{N}t\zeta\qquad
        \beta_t=\frac{\sqrt{N}t\zeta}{\alpha_t}
        \bigl(1+\frac{\norm{K}}{\alpha_t}\bigr)$$
    and choose $t$ and the rank parameter $r$ such that
    \begin{equation}\label{eq:t-admit2}
        t>1\qquad\alpha_t>0\qquad\beta_t<q_{\min}.
    \end{equation}
    Then with probability at least $1-1/t$
    \begin{equation*}
    \begin{split}
        \norm{P-\tilde P}_*<(1+\norm{\tilde G})\Bigl[
            &\frac{1}{\sqrt{q_{\min}-\beta_t}\alpha_t}
            +\frac{\sqrt{N}\tr(K)}{\sqrt{q_{\min}-\beta_t}\alpha_t^2}\\
            &+\frac{\sqrt{N}\tr(K)}{2(q_{\min}-\beta_t)^{3/2}\alpha_t^2}
            \bigl(1+\frac{\norm{K}}{\alpha_t}\bigr)\Bigr]t\zeta
    \end{split}
    \end{equation*}
    with $\tilde d_{\min}>\alpha_t>0$
    and $\tilde q_{\min}>q_{\min}-\beta_t>0$.
\end{corollary}

\begin{remark}\label{rem:bistoch-error2}
\ 
\begin{enumerate}
\item On the high probability event used in
    Corollary \ref{cor:bistoch-error2}, the positivity of
    $\tilde d_i$ and $\tilde q_i$ is ensured by the positivity
    of $d_i$ and $q_i$ and the condition \eqref{eq:t-admit2}.
\item The admissibility condition \eqref{eq:t-admit2} for $t$
    is used to ensure that $\alpha_t>0$ and $\beta_t<q_{\min}$,
    which in turn imply that $\tilde d_{\min}>0$
    and $\tilde q_{\min}>0$.
    As in Remark \ref{rem:symm-error2}, 
    the feasibility of condition \eqref{eq:t-admit2}
    can be ensured by picking a sufficiently large rank parameter $r$.
    This is because the value of $r$ can be used to control the
    expected error bound $\zeta$, which in turn determines the size
    of the admissible window of values of $t$.
\item The choice of $t$ trades the tightness of the bound against the
    failure probability $1/t$.
    As $t\to 1^+$ the error bound becomes tighter but the failure
    probability approaches 1; on the other hand, as $t$ is increased
    the failure probability becomes smaller but the error bound
    becomes looser.
\item In addition to $r$, the expected error bound $\zeta$ can also
    be controlled by the kernel function $k$ used to define $K$,
    as outlined in Remark \ref{rem:symm-error2}.
\end{enumerate}
\end{remark}

\section{Application}\label{sec:application}
We apply our algorithm for the approximate EVD of bistochastic
normalized kernel matrices to the extraction of patterns from
spatiotemporal dynamics.
As the dynamical model we consider the Kuramoto--Sivashinsky (KS)
equation
\begin{equation}\label{eq:ks}
    \partial_tu=-u\partial_xu-\partial_x^2u-\partial_x^4u\qquad
        t\geq0,\quad x\in\Xbb
\end{equation}
with periodic boundary conditions on the spatial domain
$\Xbb=[-L/2,L/2]$, $L>0$.
In the above, $u$ denotes the real valued state variable
$u(t,\cdot)\in\Ubb$, $t\geq 0$,
with state space $\Ubb\subset L^2(\Xbb,\nu)$
and $\nu$ the Lebesgue measure.

The KS equation is a dissipative partial differential equation
generating spatiotemporal chaotic dynamics.
The bifurcation parameter controlling the complexity of the dynamics
is the domain length $L$, with dynamics ranging from steady solutions
and traveling waves for low values of $L$, all the way to
spatiotemporal chaos for larger values
\cite{Kevrekidis1990,Papageorgiou1991,Cvitanovic2010}.

In addition to its rich dynamics and well understood bifurcation
diagram, the KS problem \eqref{eq:ks}
has additional desirable properties that make it an excellent testbed
for pattern extraction methods.
First, its solutions satisfy dynamical symmetries;
more specifically, the KS problem \eqref{eq:ks}
is equivariant under spatial translations
$u(t,x)\mapsto u(t,x+y)$ for all $y\in\Rbb$
and anti-reflection $u(t,x)\mapsto -u(t,-x)$.
Second, it has a global attractor of finite dimension;
namely, a finite dimensional subset of $\Ubb$
which is forward invariant and attracts almost all
initial conditions $u(0,\cdot)\in\Ubb$
\cite{Temam1997,RobinsonJC2001}.

The dynamics generated by the KS problem \eqref{eq:ks}
is given by the flow map $\Phi^t\colon\Ubb\to\Ubb$,
$\Phi^t(u(t_0,\cdot))=u(t_0+t,\cdot)$,
with continuous time variable $t\in\Rbb$,
state space $\Ubb$ and invariant probability measure
$\mu$ with compact support.
In what follows, we also employ the discrete time flow map
$\Phi^n=\Phi^{n\Delta t}$ with sampling timestep
$\Delta t>0$ and $n\in\Zbb$.
The observables of the dynamics are members of the
Hilbert space of real valued functions $L^2(\Ubb,\mu)$.

\subsection{Spatiotemporal pattern extraction}
The problem of identifying spatiotemporal patterns of a dynamical
system has traditionally been formulated as an eigendecomposition
problem for a kernel integral operator acting on the space of
observables $L^2(\Ubb,\mu)$.
The most popular method is arguably
proper orthogonal decomposition (POD),
where the kernel integral operator is formed using a two-point
correlation kernel \cite{Aubry1991,Berkooz1993,Hinze2005}.
In this work we employ an alternative but related approach
called vector valued spectral analysis (VSA)
\cite{Giannakis2019vsa}.

The VSA method employs the product state space
$\Omega=\Ubb\times\Xbb$
and associated real Hilbert space
$H=L^2(\Omega,\sigma)$
with product measure $\sigma=\mu\times\nu$.
Every function $f\in H$ represents a spatiotemporal pattern
of the dynamics, with a temporal dependence through $u\in\Ubb$
and a spatial dependence through $x\in\Xbb$.
More specifically, for every $u\in\Ubb$,
$f(u,\cdot)\in L^2(\Xbb,\nu)$
denotes a function on the spatial domain $\Xbb$,
with $f(u,x)\in\Rbb$ its pointwise value at $x\in\Xbb$.
The map $t\mapsto f(\Phi^t(u),\cdot)$
represents the temporal evolution of a pattern $f\in H$
by the dynamics $\Phi^t$ for an initial state $u\in\Ubb$.

The desired spatiotemporal patterns are given by the eigenfunctions
of a kernel integral operator
$\mathcal{K}\colon H\to H$
\begin{equation}\label{eq:integral-op}
    \mathcal{K}f(\omega)=\int_\Omega\kappa(\omega,\omega')
        f(\omega')d\sigma(\omega')
\end{equation}
with product state $\omega=(u,x)\in\Omega$
and continuous, bounded and positive semidefinite kernel function
$\kappa\colon\Omega\times\Omega\to\Rbb$.
The integral operator $\mathcal{K}$ is compact and self-adjoint;
as a result, its eigenfunctions can be chosen to form an orthonormal
basis of $H$, with eigenvalues that are real, nonnegative and have zero
as their limit point.
By forming an operator acting directly on $H$ we obtain
spatiotemporal patterns that are generally not of tensor product
form, meaning that they are not expressible as the tensor product
of a pair of temporal and spatial modes.
This is in contrast to traditional approaches such as POD,
where one computes the eigenfunctions of an operator acting on the
temporal space $L^2(\Ubb,\mu)$ and forms their tensor product
with a basis for the spatial space $L^2(\Xbb,\nu)$.

\begin{figure}[t]
    \centering
    \includegraphics[width=.8\textwidth]{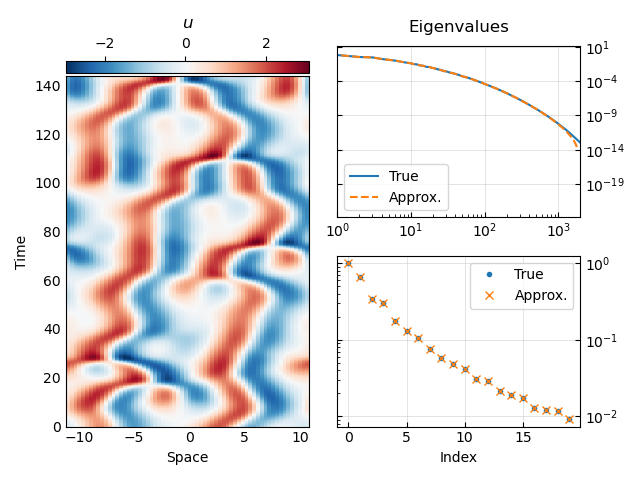}
    \caption{(Left) Space-time heatmap of the true state data
    obtained by integrating the KS problem \eqref{eq:ks}
    with $L=7\pi$ for $575$ timesteps using the parameter values
    given in Section \ref{sec:numerics}.
    (Right) Comparison of the true eigenvalues with those
    obtained using Algorithm \ref{alg:bistoch} with
    $r=2\,048$.
    The horizontal axis of the top panel is limited to the
    maximum value of $2\,000$ to facilitate the comparison.
    The bottom panel focuses on the leading 20 eigenvalues.}
    \label{fig:compare1}
\end{figure}

For our kernel $\kappa=k\circ(W\times W)$
we employ a Gaussian kernel function
$k\colon\Rbb^J\times\Rbb^J\to\Rbb$
acting by composition with a delay embedding map
$W\colon\Omega\to\Rbb^J$,
where $J\in\Nbb$ denotes the number of time delays.
Given a product state sample $\omega=(u,x)\in\Omega$,
the delay embedding map $W$ forms $J$ delays in time
at the spatial point $x\in\Xbb$
$$W((u,x))=(u(x),\Phi^{-1}(u)(x),\ldots,
    \Phi^{-(J-1)}(u)(x))$$
and the Gaussian kernel function $k$ acts by
\begin{equation}\label{eq:data-kernel}
    k(W(\omega),W(\omega'))=\exp\bigl(-\frac{1}{\varepsilon J}
        \norm{W(\omega)-W(\omega')}_2^2\bigr)
\end{equation}
where $\varepsilon>0$ is a tunable bandwidth parameter
and $\norm{\cdot}_2$ denotes the 2-norm in $\Rbb^J$.
For our numerical experiments we will use the
bistochastic normalized versions of kernel \eqref{eq:data-kernel}
and associated kernel integral operator \eqref{eq:integral-op}.

\begin{figure}[t]
    \centering
    \includegraphics[width=.9\textwidth]{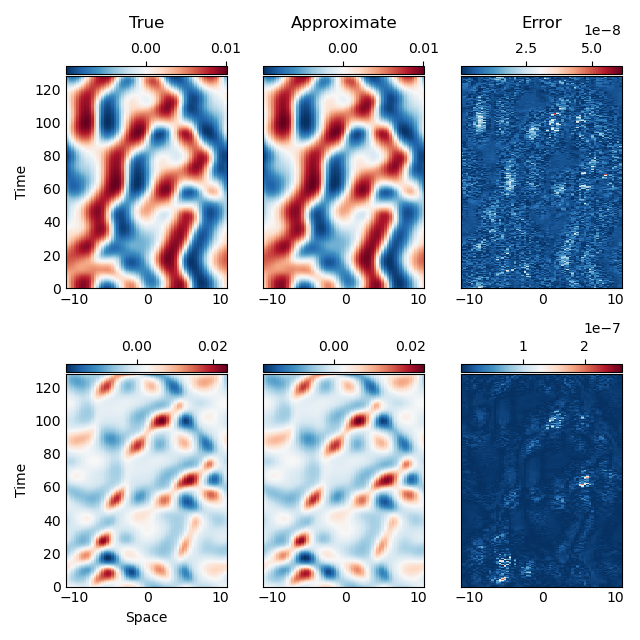}
    \caption{Comparison of the true eigenfunctions
    $\phi_1$ (top row) and $\phi_2$ (bottom row)
    with those obtained using Algorithm \ref{alg:bistoch}
    with $r=2\,048$.}
    \label{fig:compare2}
\end{figure}

Because it acts via composition with the delay embedding map $W$,
kernel $\kappa$ factors the product state space $\Omega$
into equivalence classes consisting of states with identical
dynamical behavior under $J$ delays.
As shown in \cite{Giannakis2019vsa},
this implies that the functions in the range $\ran\mathcal{K}$
of the integral operator \eqref{eq:integral-op}
are invariant under the actions of spatial symmetries of the
KS problem \eqref{eq:ks}.
To make this property precise, we consider the group of symmetries
$G$ with continuous left action $\Gamma_g$ on the spatial domain
$\Xbb$ for every $g\in G$.
Every induced action $\Gamma_{\Ubb,g}$ on $\Ubb$,
$\Gamma_{\Ubb,g}(u)=u\circ\Gamma_g^{-1}$,
represents a dynamical symmetry of the dynamics generated
by \eqref{eq:ks}.
This means that the dynamics $\Phi^t$ satisfies the
equivariance property
$$\Gamma_{\Ubb,g}\circ\Phi^t=\Phi^t\circ\Gamma_{\Ubb,g}$$
for all $g\in G$ and $t\geq 0$.
For the KS problem \eqref{eq:ks}, the induced actions
$\Gamma_{\Ubb,g}$ represent spatial translations of the state
variable $u(t,\cdot)$.
Note that the anti-reflection symmetry of \eqref{eq:ks}
cannot be induced by an action $\Gamma_g$ acting on the spatial
domain $\Xbb$.

For our choice of kernel function $\kappa$, every function
$f\in\ran\mathcal{K}$ satisfies the analogous invariance property
\begin{equation}\label{eq:G-invariance}
    f\circ\Gamma_{\Omega,g}=f
\end{equation}
with induced action
$\Gamma_{\Omega,g}=\Gamma_{\Ubb,g}\times\Gamma_g$
for all $g\in G$.
Thanks to the invariance property \eqref{eq:G-invariance},
the eigenfunctions of the kernel integral operator
\eqref{eq:integral-op}
are constant on product states in $\Omega$ that differ only
by a symmetry transformation.
As a result, each such function can generally represent a more
complex spatiotemporal pattern than when symmetry invariance is
not ensured \cite{Giannakis2019vsa}.

\begin{figure}[t]
    \centering
    \includegraphics[width=.9\textwidth]{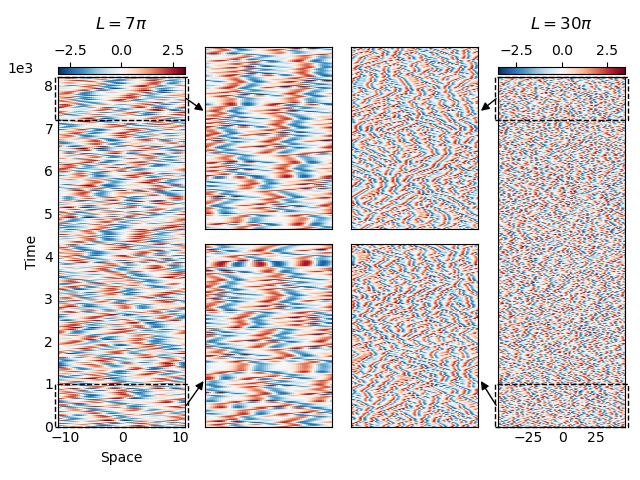}
    \caption{(Two leftmost columns) Space-time heatmaps of the
    true state data obtained by integrating the KS problem
    \eqref{eq:ks} with $L=7\pi$ for $32\,768$ timesteps
    using the parameter values given in
    Section \ref{sec:numerics}.
    (Two rightmost columns) Space-time heatmaps for the
    corresponding results with $L=30\pi$.}
    \label{fig:state}
\end{figure}

\subsection{Numerical experiments}\label{sec:numerics}
For our numerical experiments we consider the KS problem
\eqref{eq:ks} with domain lengths $L=7\pi$ ($\approx 22$)
and $L=30\pi$ ($\approx 94$),
both of which lead to chaotic dynamics
\cite{Kevrekidis1990,Papageorgiou1991}.
For the $L=7\pi$ case, we perform a Fourier spatial
discretization using $M=64$ Fourier modes and 3/2 dealiasing
of the pseudospectral treatment of the quadratic nonlinearity.
For the temporal discretization we employ the exponential
time differencing 4-stage Runge--Kutta method
\cite{Cox2002,Kassam2005}
with timestep $\Delta t=0.25$.
Our initial condition is formed by setting the leading four
Fourier coefficients to $0.6$ and the rest to zero.
Our simulations with $L=30\pi$ follow the exact same approach,
except that $M=128$ Fourier modes are used for the spatial
discretization.

We begin by considering a training dataset of relatively small
size, which will allow us to compare the patterns extracted by
Algorithm \ref{alg:bistoch}
with the true patterns extracted by computing the true
EVD of the bistochastic version of the discretized operator
\eqref{eq:integral-op}.
This is performed for the domain length $L=7\pi$.
After integrating the discretized dynamics for $10\,000$ timesteps
($2\,500$ time units), we collect one sample every timestep
for a total of $575$ time samples.
A space-time plot of the obtained solution is shown in the
left panel of Figure \ref{fig:compare1}.
Using $J=64$ delays, the resulting dataset consists of $N=512$
time samples in delay embedded form, bringing the total number of
product state samples to $NM=32\,768$.\footnote{Throughout this
section we use $NM$ to denote the total number of samples,
consisting of $N$ temporal and $M$ spatial samples;
this deviates from our notation in other sections of the paper,
where $N$ denotes the total number of samples.}
The bandwidth parameter $\varepsilon$ of kernel \eqref{eq:data-kernel}
is calibrated using the scaling-based algorithm developed in
\cite{Coifman2008}.
All presented results for $L=7\pi$ use the bandwidth value
$\varepsilon=0.8$.

\begin{figure}[t]
    \centering
    \includegraphics[width=.9\textwidth]{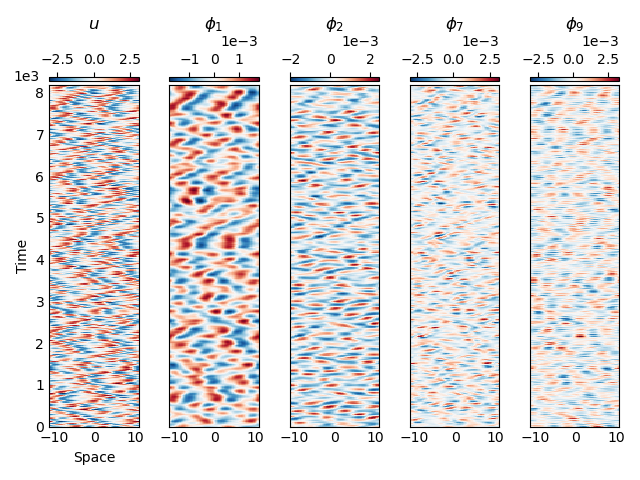}
    \caption{Space-time heatmaps of the true state data $u$
    and select eigenfunctions $\phi_i$ for the $L=7\pi$
    simulation data.}
    \label{fig:eigfun1a}
\end{figure}

We employ Algorithm \ref{alg:bistoch} to approximate
the EVD of the kernel matrix corresponding to the bistochastic
normalized version of kernel \eqref{eq:data-kernel}.
For that we perform a partial Cholesky factorization using ARPC-LM
with rank parameter $r=2\,048$ and block size $b=64$,
which yields an approximation with realized rank $1\,789$
and relative trace norm error $0.01\%$.
As explained in Appendix \ref{app:rpc},
ARPC-LM samples pivots in blocks and subjects them to an
additional accept/reject check.
This means that the final number of accepted pivots (the realized
rank) can be lower than the proposed rank $r$.
This computation was carried out using one Nvidia A100 40 GiB GPU.
The computed eigenvalues and eigenfunctions are compared with the true 
ones obtained by directly computing the EVD of the bistochastic
kernel matrix.
The true EVD was computed on an AMD EPYC 7713 CPU.

The two right panels of Figure \ref{fig:compare1} compare the
true and approximate eigenvalues, showing the close agreement
obtained up to the realized rank of the employed approximation.
As expected, the accuracy of the approximation starts to deteriorate
close to the tail of the captured eigenvalues.
The onset of this deterioration is controlled by the realized rank
of the approximation.
Figure \ref{fig:compare2}
compares the second ($\phi_1$) and third ($\phi_2$)
eigenfunctions obtained using
Algorithm \ref{alg:bistoch}
with the corresponding true ones.
Since we are using an ergodic bistochastic kernel integral
operator, the leading eigenfunction $\phi_0$ is a constant
function in both the approximate and true results,
which is why we do not include it in our comparisons.
Because our training data consists of only one dynamical
trajectory, we can make the identification $t\mapsto u(t,\cdot)$,
which allows us to plot each $\phi_i(\cdot,x)$
as a function of time.
The results demonstrate the close pointwise agreement between
the true and approximate eigenfunctions shown.
The agreement persists to higher eigenfunctions in a similar way
to the agreement between the associated eigenvalues shown in
Figure \ref{fig:compare1}.

\begin{figure}[t]
    \centering
    \includegraphics[width=.9\textwidth]{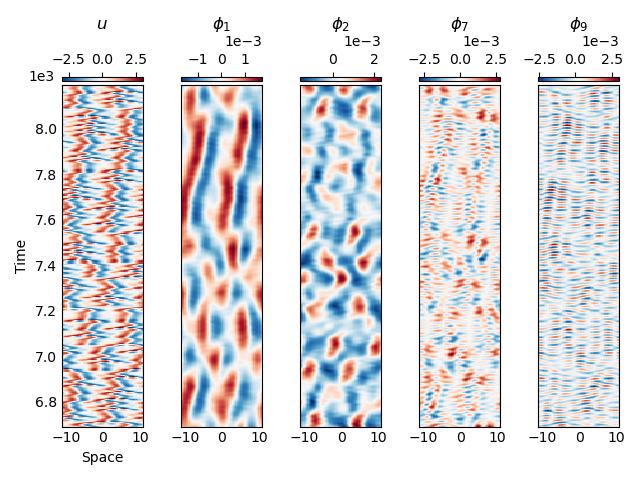}
    \caption{As in Figure \ref{fig:eigfun1a},
    but restricted to the last $1\,500$ time units.}
    \label{fig:eigfun1b}
\end{figure}

To facilitate the comparison of the approximate EVD with the
true one, we have so far restricted ourselves to a dataset of
small size with $L=7\pi$ only.
We now consider both domain lengths $L=7\pi$ and $L=30\pi$
and increase the dataset size to $N=32\,768$ time samples
and the number of delays to $J=512$, keeping the timestep
and spatial resolution unchanged.
As a reminder, the $L=7\pi$ results use $M=64$ Fourier modes
for the spatial discretization; those with $L=30\pi$ use
$M=128$ Fourier modes.
The true state data obtained by integrating the KS problem
is shown in Figure \ref{fig:state}.
Both domain lengths yield spatiotemporal chaos, with the larger
spatial domain leading to finer spatial features in the dynamics.
All computations that follow were performed using four
Nvidia A100 40 GiB GPUs.

We begin by considering the $L=7\pi$ results.
The collected state samples produce a product space of
dimension $NM=2\,097\,152$,
making the direct computation of the true EVD unfeasible on
most devices.
To compute the approximate EVD, we use ARPC-LM with
block size $b=512$ and rank parameter $r=32\,768$,
yielding an approximation with realized rank $32\,057$
and relative trace norm error $7.51\%$.
We then employ Algorithm \ref{alg:bistoch}
to extract the leading $512$ eigenfunctions of the
bistochastic normalized kernel matrix.

Figures \ref{fig:eigfun1a} and \ref{fig:eigfun1b}
show a collection of the obtained eigenfunctions,
using the same identification $t\mapsto u(t,\cdot)$
employed earlier.
Figure \ref{fig:eigfun1a} shows the eigenfunctions
$\phi_1$, $\phi_2$, $\phi_7$ and $\phi_9$ along with the
true state data.
Figure \ref{fig:eigfun1b} presents the same results,
focusing on the values of the true state data and the
selected eigenfunctions for the last $1\,500$ time units
of the full time horizon.
Since the leading eigenfunction $\phi_0$ (not shown)
is constant, all higher eigenfunctions are zero-mean functions
that encode the variation of the dynamics relative to the mean
captured in $\phi_0$.
The eigenfunctions presented in Figures
\ref{fig:eigfun1a}--\ref{fig:eigfun1b}
encode patterns resembling perturbed traveling waves,
with higher eigenfunctions generally leading to finer spatial
variation and faster temporal evolution.

\begin{figure}[t]
    \centering
    \includegraphics[width=.9\textwidth]{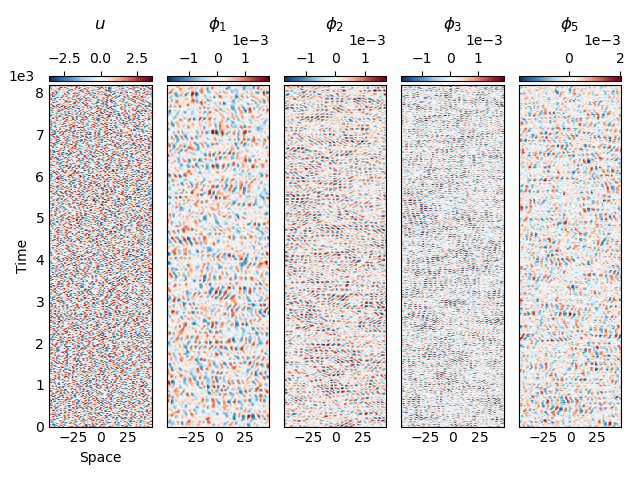}
    \caption{Space-time heatmaps of the true state data $u$
    and select eigenfunctions $\phi_i$ for the $L=30\pi$
    simulation data.}
    \label{fig:eigfun2a}
\end{figure}

The ability to encode complex spatiotemporal patterns of the
underlying dynamics in individual eigenfunctions is a result of
the invariance property \eqref{eq:G-invariance}
satisfied by the VSA kernel eigenfunctions.
We refer the reader to \cite{Giannakis2019vsa} for more details
on this point, which includes a comparison of VSA eigenfunctions
with those obtained using POD, showing that the POD patterns represent
the evolution of pure Fourier modes, which individually have little
value in explaining important patterns of the dynamics.
The patterns presented in the present work are in qualitative
agreement with those obtained for $L=22$ in \cite{Giannakis2019vsa},
where similar KS simulation parameters are used,
but the normalization of the employed kernel is different.

Next we consider the case of the larger spatial domain $L=30\pi$,
leading to finer spatial variation in the produced chaotic dynamics
as shown in Figure \ref{fig:state}.
Due to the increased number $M=128$ of Fourier modes used for the
spatial discretization, this dataset leads to a product space
of dimension $NM=4\,194\,304$.
We again employ ARPC-LM and Algorithm \ref{alg:bistoch}
to extract the leading $512$ eigenfunctions, using the kernel
bandwidth value $\varepsilon=1.2$.
For ARPC-LM we use block size $b=512$ and rank $r=32\,768$,
producing an approximation with realized rank $32\,215$
and relative trace norm error $14.4\%$.

\begin{figure}[t]
    \centering
    \includegraphics[width=.9\textwidth]{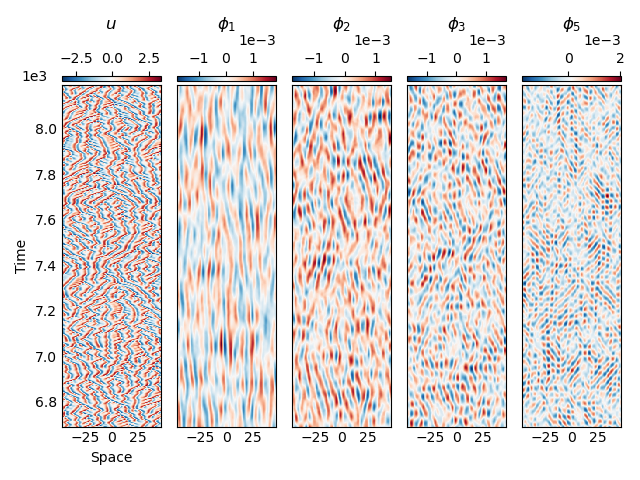}
    \caption{As in Figure \ref{fig:eigfun2a},
    but restricted to the last $1\,500$ time units.}
    \label{fig:eigfun2b}
\end{figure}

Figures \ref{fig:eigfun2a} and \ref{fig:eigfun2b}
present the true state data and eigenfunctions $\phi_1$, $\phi_2$,
$\phi_3$ and $\phi_5$, with Figure \ref{fig:eigfun2b}
focusing on their values for the last $1\,500$ time units.
As in the $L=7\pi$ case, the extracted eigenfunctions resemble
perturbed traveling waves that become increasingly fine in their
spatiotemporal variation as the eigenfunction index is increased.
Similar features are captured by the higher eigenfunctions not
shown here.
The extracted patterns are in qualitative agreement with those
obtained for $L=94$ in \cite{Giannakis2019vsa}
using a smaller dataset and different kernel normalization.

\begin{figure}[t]
    \centering
    \includegraphics[width=.9\textwidth]{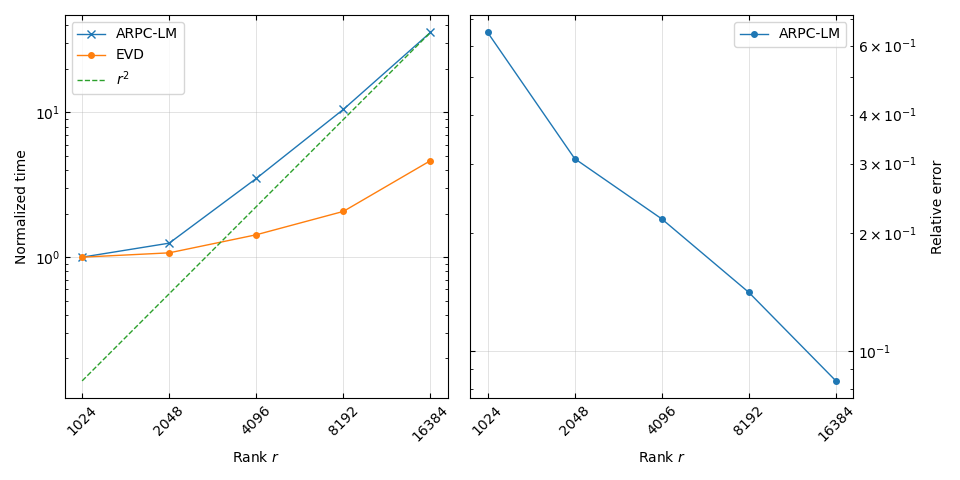}
    \caption{(Left) Scaling of the normalized wall-clock time for
    ARPC-LM and Algorithm \ref{alg:bistoch} (EVD, excluding ARPC-LM)
    with the rank parameter $r$.
    The two data curves are normalized independently by their own
    wall-clock time at $r=1\,024$, which is why both of them are
    equal to 1 at that $r$ value.
    (Right) Variation of the relative trace norm error of ARPC-LM
    with $r$.}
    \label{fig:scaling}
\end{figure}

Finally, we perform a preliminary investigation of the scaling
of the wall-clock execution time and relative error for
ARPC-LM and Algorithm \ref{alg:bistoch}
with the rank parameter $r$.
For that we use an $L=7\pi$ dataset consisting of
$N=32\,768$ delay embedded samples and $M=64$ spatial gridpoints,
with $J=64$ delays.
The results are shown in Figure \ref{fig:scaling}
and were produced using two Nvidia A100 40 GiB GPUs
and constant ARPC-LM block size $b=256$
for all tested rank parameters $r$.
In the left figure panel we see that the execution time of ARPC-LM
is scaling approximately quadratically for the largest $r$ values
in our tested range.
The execution time of Algorithm \ref{alg:bistoch}
(EVD, excluding ARPC-LM) is increasing sub-quadratically
for most of the tested $r$ values, but also seems to approach
the expected quadratic scaling toward the end of the tested range.
As $r$ is increased beyond our current range, we expect
the EVD curve to approach a quadratic scaling as well.
The sub-quadratic scaling is likely present because the computation
is limited by the inter-device communication rate at low $r$ values,
not by the compute rate.
The maximum $r$ value tested in this preliminary investigation
is set by the available memory of the employed GPU devices.
The right figure panel shows the variation of the relative trace
norm error obtained by ARPC-LM for the tested range of $r$.

\begin{remark}\label{rem:results}
\ 
\begin{enumerate}
\item The choice of kernel function used in defining the kernel
    integral operator \eqref{eq:integral-op}
    plays an important role in determining its rank,
    and ultimately the maximum rank that can be obtained by the
    linear span of the numerically computed eigenfunctions.
    Our use of a Gaussian kernel means that the rank of the integral
    operator is infinite, regardless of the rank of the employed data
    matrix \cite{Schoelkopf2001,Giannakis2019vsa}.
    In contrast, the state correlation kernels traditionally used
    in POD do not share that property, meaning that the rank of the
    associated integral operator is bounded by the rank of the data
    matrix.
    Although the rank obtained by POD is sufficient to represent
    the state data to any given level of accuracy in an $L^2$ norm,
    it might not be sufficient to represent arbitrary observables.
    As a result, our choice of a Gaussian kernel is motivated by
    applications where representing observables by the
    computed eigenfunctions is important
    \cite{Freeman2023,Freeman2024,Vales2026qmcl}.
\item For our numerical results we use the Gaussian kernel
    \eqref{eq:data-kernel} with fixed bandwidth $\varepsilon$.
    The methods considered in this work can also be applied to
    kernels of variable bandwidth, where the employed bandwidth
    depends on the arguments of the kernel function
    \cite{Berry2016a}.
\end{enumerate}
\end{remark}

\section{Conclusion}\label{sec:conclusion}
We developed an algorithm for the approximate computation
of the eigenvalue decomposition of normalized kernel matrices,
focusing on a symmetric and a bistochastic normalization.
The proposed algorithm employs a pivoted partial Cholesky
algorithm to construct a low rank approximation of the original
kernel matrix and compute the approximate eigenvalue decomposition
of its normalization, relying on a limited number of kernel evaluations.
We applied the developed algorithm to the kernel based extraction
of spatiotemporal patterns from chaotic dynamics, and investigated
its accuracy and scalability.
In addition to numerical results, we derived trace norm error bounds
for the proposed approximations.

Next steps in this line of research involve using our
proposed algorithm to enable the application of normalized
kernel methods to large datasets for tasks such as spatiotemporal
pattern extraction, model reduction and dynamical closure.

\appendix
\section{Randomly pivoted Cholesky}\label{app:rpc}
In this section we provide a brief overview of the randomly
pivoted Cholesky (RPC) algorithm and the accelerated variant
\cite{ChenY2025,Epperly2025} used in the present work.
We use zero-based indexing notation to denote entries of matrices.

Following the material in Section \ref{sec:intro},
we consider an $N\times N$ positive semidefinite kernel matrix
$K$ resulting from the evaluation of the kernel function $k$
on a dataset $\{u_i\}_{i=0}^{N-1}$.
Given a rank parameter $r<N$, RPC chooses $r$ column indices
(pivots) of $K$ and uses the associated columns to compute the
partial Cholesky factorization
$$K\approx\tilde{K}=FF^\top$$
with $N\times r$ partial Cholesky factor $F$.
The factor matrix $F$ is computed iteratively.
In each iteration $0\leq i<r$, RPC selects a pivot
$0\leq s_i<N$,
evaluates the associated kernel matrix column $K(:,s_i)$
and computes the new column $F(:,i)$.

To sample the pivots and measure the approximation error in the
trace norm in each iteration $i$, the algorithm keeps track of
the diagonal of the residual matrix
$$\rho^{(i)}=\diag(K-\tilde{K}^{(i)})$$
where $\tilde{K}^{(i)}$ denotes the approximate kernel matrix
in iteration $i$.
Using the vector $\rho^{(i)}$, the pivot $s_i$ is sampled
according to the discrete probability distribution
$$\Pbb\{s_i=j\}=\frac{\rho^{(i)}(j)}{\sum_{k=0}^{N-1}\rho^{(i)}(k)}
    \qquad j\in\{0,\ldots,N-1\}.$$
Namely, the algorithm uses the diagonal entries of the residual
matrix to inform its pivot sampling in each iteration.

The use of the diagonal entries is motivated by the fact that the
off-diagonal entries of a positive semidefinite matrix are bounded
in absolute value by those on the diagonal.
As a result, a large diagonal entry $\rho^{(i)}(j)$
means that entries of large absolute value may be present
in column $(K-\tilde K^{(i)})(:,j)$,
suggesting that the column index $j$ may be a good next pivot choice.
This randomized pivot selection can be seen as a compromise between
two extremes:
(1) the deterministic greedy case, where in each iteration one
chooses the pivot $s_i$ corresponding to the largest entry on the
diagonal $\rho^{(i)}$;
(2) the uniformly random case, where one draws a pivot $s_i$
uniformly at random, without using any information about the
diagonal entries.

\begin{algorithm}[t]
    \caption{Accelerated RPC, low memory (ARPC-LM) \cite{Epperly2025}}
\label{alg:arpclm}
\begin{algorithmic}
    \STATE Input: $N\times N$ kernel matrix $K$,
        maximum rank $r<N$, block size $b$
    \STATE Output: accepted pivots $\Ssf\subseteq\{s_0,\ldots,s_{r-1}\}$,
        $|\Ssf|\times |\Ssf|$ inverse factor $L^{-1}$
    \STATE
    \STATE Initialize $L^{-1}\gets 0$, $\Ssf\gets\emptyset$,
        $\rho\gets\texttt{diag}(K)$
    \STATE Ensure $N$ divisible by $b$, $r$ divisible by $b$
    \FOR{$i=0$ \TO $r/b-1$}
        \STATE $s_{ib},\ldots,s_{(i+1)b-1}\sim\rho/\texttt{sum}(\rho)$
            \hfill iid with replacement
        \STATE $\Ssf_i'\gets\{s_{ib},\ldots,s_{(i+1)b-1}\}$
            \hfill proposed pivots
        \STATE Evaluate submatrix $K(\Ssf_i',\Ssf\cup\Ssf_i')$
        \STATE $H\gets K(\Ssf_i',\Ssf_i')-[L^{-1}K(\Ssf_i',\Ssf)^\top]^\top
            [L^{-1}K(\Ssf_i',\Ssf)^\top]$
        \STATE $\Ssf_i, L_i^{-1}\gets \texttt{reject}(H,\Ssf_i')$\hfill
            accept/reject check (Alg.\@ \ref{alg:reject})
        \STATE Extract $K(\Ssf,\Ssf_i)$ from
            $K(\Ssf_i',\Ssf\cup\Ssf_i')^\top$\hfill
            available since $\Ssf_i\subseteq \Ssf_i'$
        \STATE $M\gets L^{-\top}L^{-1}K(\Ssf,\Ssf_i)$
            \hfill precompute; independent of $\Rsf$
        \FOR{$j=0$ \TO $N/b-1$}
            \STATE $\Rsf\gets\{jb,\ldots,(j+1)b-1\}$
                \hfill block of $b$ rows
            \STATE Evaluate submatrix $K(\Rsf,\Ssf\cup\Ssf_i)$
            \STATE $G\gets (K(\Rsf,\Ssf_i)-K(\Rsf,\Ssf)M)L_i^{-\top}$
            \STATE $\rho(\Rsf)\gets\rho(\Rsf)-\texttt{rownorms}(G)^2$\hfill
                update residual
        \ENDFOR
        \STATE $L^{-1}\gets\begin{bmatrix}L^{-1} & 0\\
            -L_i^{-1}M^\top & L_i^{-1}\end{bmatrix}$\hfill
            update inverse Cholesky factor
        \STATE $\rho\gets\texttt{max}(\rho,0)$\hfill ensure nonnegativity
        \STATE $\Ssf\gets\Ssf\cup\Ssf_i$\hfill update pivots
    \ENDFOR
\end{algorithmic}
\end{algorithm}

In this work we use the low memory, accelerated version (ARPC)
of the algorithm, which samples pivots in blocks instead of
one by one \cite{Epperly2025}.
Given a block size $b\in\Nbb$, ARPC samples $b$ pivots
in each round based on the probability distribution described above.
The sampled pivots undergo an additional check, where the algorithm
uses a randomized procedure to decide whether to accept or reject
the sampled pivots.
Only the accepted pivots are used for building the factor matrix $F$.
The accept/reject check of ARPC ensures that, in each round, pivots are
accepted with the same probability as in RPC, where only one pivot
is drawn in each round.

\begin{algorithm}[t]
    \caption{ARPC-LM accept/reject check \cite{Epperly2025}}
\label{alg:reject}
\begin{algorithmic}
    \STATE Input: $b\times b$ matrix $H$,
        proposed pivots $\Ssf'=\{s_0',\,\ldots,\,s_{b-1}'\}$
    \STATE Output: accepted pivots $\Ssf\subseteq\Ssf'$,
        $|\Ssf|\times |\Ssf|$ inverse factor $L^{-1}$
    \STATE
    \STATE $L\gets 0$, $\Tsf\gets\emptyset$
    \STATE $\rho\gets\texttt{diag}(H)$
    \FOR{$i=0$ \TO $b-1$}
        \STATE $z\gets\texttt{randu}(0,\,1)$
            \hfill uniformly random in $[0,1)$
        \IF{$z\cdot\rho(i)<H(i,i)$}
            \STATE $\Tsf\gets\Tsf\cup\{i\}$\hfill accept pivot
            \STATE $L(i:b,i)\gets H(i:b,i)/\sqrt{H(i,i)}$
                \hfill update Cholesky factor
            \STATE $H(i+1:b,i+1:b)\gets H(i+1:b,i+1:b)
                -L(i+1:b,i)L(i+1:b,i)^\top$
        \ENDIF
    \ENDFOR
    \STATE $L^{-1}\gets L(\Tsf,\Tsf)^{-1}$\hfill
        extract accepted entries \& invert
    \STATE $\Ssf\gets\{s_i':i\in\Tsf\}$\hfill extract accepted pivots
\end{algorithmic}
\end{algorithm}

The low memory variant (ARPC-LM) of ARPC does not store
the full $N\times r$ factor matrix $F$ in memory.
This is done to remove the dependence of the required memory
on the full dataset size $N$, which makes the application
of ARPC prohibitively memory expensive for large datasets.
Instead, ARPC-LM relies on the $r\times r$
Cholesky factor $L$ introduced in \eqref{eq:nystrom-factored}.
As explained in Section \ref{sec:cholesky},
the pivoted partial Cholesky algorithm can be used to compute
the Nystr\"om approximation $\tilde K$ in the factored form
$$\tilde K=FF^\top=\bigl[K(:,\Ssf)L^{-\top}\bigr]
    \bigl[K(:,\Ssf)L^{-\top}\bigr]^\top.$$
Namely, the partial Cholesky factor $F=K(:,\Ssf)L^{-\top}$
can always be reconstructed from the Cholesky factor $L$
and the data sampled by the pivot set $\Ssf$.
Therefore, ARPC-LM stores the $r\times r$ matrix $L$
and recomputes the entries of the kernel matrix $K$
as required.
This removes the dependence of the required memory on
the dataset size $N$, while increasing the total number
of kernel evaluations; nevertheless, the overall
computational cost of the method remains $O(Nr^2)$.
In addition, the full and low memory versions of ARPC are
theoretically equivalent, which means that they propose
and accept pivots with the same probability distributions
in each round \cite{Epperly2025}.

Our implementation of ARPC-LM is presented
in pseudocode form in Algorithm \ref{alg:arpclm},
with the accept/reject check in Algorithm \ref{alg:reject}.
Algorithm \ref{alg:arpclm} is a modified version of the
original algorithm \cite[Algorithm 2.3]{Epperly2025};
the main difference is that the present version works with the
inverse Cholesky factor $L^{-1}$ directly, instead of $L$.
This is motivated by the fact that all matrix operations in
Algorithm \ref{alg:arpclm} and the downstream EVD approximation
employ $L^{-1}$, never $L$.

ARPC-LM takes as inputs the kernel matrix $K$, rank parameter $r$
and block size $b$, and outputs the inverse Cholesky factor $L^{-1}$
and accepted pivots $\Ssf\subseteq\{s_0,\ldots,s_{r-1}\}$.
Note that because pivots are drawn in blocks, which then undergo the
additional accept/reject check, the algorithm does not guarantee that
the square matrix $L^{-1}$ will have exactly $r$ rows/columns,
only that its number of rows/columns is bounded from above by $r$.
The same is true for the total number of accepted pivots $|\Ssf|$.
In addition, the algorithm does not evaluate all entries of $K$
and does not need to store the full kernel matrix in memory.

As with the implementation of Algorithm \ref{alg:bistoch},
to generate our results in Section \ref{sec:application}
we use a distributed GPU implementation of
Algorithm \ref{alg:arpclm}.
In the employed code, each GPU holds a block of columns of the
$r\times r$ inverse Cholesky factor $L^{-1}$,
enabling the scaling of the algorithm to large
datasets and approximation ranks where even the lower-dimensional
$r\times r$ factor is too large to fit in the memory of one GPU.
The inner loop of Algorithm \ref{alg:arpclm}
used for updating all $N$ entries of the residual vector $\rho$
is split among the available GPUs, with each one being responsible
for updating its share of the $N$ rows.
All $N$-dependent computations such as evaluating large blocks of
kernel entries or updating the residual vector $\rho$
are performed in single precision (32-bit) floating point arithmetic.
The rest, including the accept/reject check of
Algorithm \ref{alg:reject}
and the iterative formation of $L^{-1}$,
are performed in double precision (64-bit).
Additional details can be found in the accompanying
code\footnote{\url{https://github.com/cval26/kernel_evd}}.

To simplify the presentation, in Section \ref{sec:approx}
and Algorithm \ref{alg:bistoch} we assume that the inverse
Cholesky factor $L^{-1}$ is exactly of size $r\times r$,
which also implies that the partial Cholesky factor
$F=K(:,\Ssf)L^{-\top}$ has exactly $r$ columns,
instead of at most $r$ columns.
This simplification does not affect any of the presented results
or analysis, as any columns of $F$ past the number of accepted
pivots can be set to zero without affecting the downstream
computations.
The numerical implementation of Algorithms
\ref{alg:bistoch} and \ref{alg:arpclm}
does not rely on that assumption.

\section{Proofs}\label{app:proof}
In this section we provide proofs for the results presented in
Section \ref{sec:error}.
Assumptions \ref{asm:kernel} and \ref{asm:psd}
are in effect throughout this section.
We use $\norm{\cdot}$ for the operator norm
and $\norm{\cdot}_*$ for the trace norm.

\subsection{Symmetric normalization}\label{app:proof-symm}
We prove the results of Section \ref{sec:error-symm}
related to the normalized matrices $L$ and $\tilde L$. 
We begin by proving the following lemma bounding the difference
between the original and approximate normalization constants.
\begin{lemma}\label{lem:dbound}
    We have
    $$\sum_{i=0}^{N-1}(d_i-\tilde d_i)^2\leq N[\tr(K-\tilde K)]^2,
        \qquad\max_{0\leq i<N}\abs{d_i-\tilde d_i}
        \leq\sqrt{N}\tr(K-\tilde K).$$
\end{lemma}
\begin{proof}
    The matrix $K-\tilde K$ is positive semidefinite, so its eigenvalues
    are nonnegative and bounded above by $\tr(K-\tilde K)$,
    which implies that $\norm{K-\tilde K}\leq\tr(K-\tilde K)$.
    As a result,
    $$\sum_{i=0}^{N-1}(d_i-\tilde d_i)^2=\norm{(K-\tilde K)1_N}_2^2
        \leq\norm{K-\tilde K}^2\norm{1_N}_2^2
        \leq N[\tr(K-\tilde K)]^2$$
    which also yields
    $\max_{0\leq i<N}\abs{d_i-\tilde d_i}\leq\sqrt{N}\tr(K-\tilde K)$.
\end{proof}

\begin{theorem2}[\ref{thm:symm-error1}]
    Under Assumption \ref{asm:dpos}, the normalized matrices
    $L$ and $\tilde L$ satisfy
    $$\lVert L-\tilde L\rVert_*\leq\Bigl[\frac{1}{d_{\min}}
        +\frac{\sqrt{N}}{\delta^2}\tr(K)
        +\frac{N}{4\delta^3}\tr(K)\tr(K-\tilde K)\Bigr]
        \tr(K-\tilde K).$$
\end{theorem2}
\begin{proof}
    Write $L-\tilde L=A+B$ with terms
    $$A=D^{-1/2}(K-\tilde K)D^{-1/2}\qquad
        B=D^{-1/2}\tilde KD^{-1/2}
        -\tilde D^{-1/2}\tilde K\tilde D^{-1/2}.$$
    For matrix $A$, we write
    $$\lVert A\rVert_*=\tr(A)=\tr[D^{-1}(K-\tilde K)]
        \leq\lVert D^{-1}\rVert\tr(K-\tilde K)
        \leq\frac{1}{d_{\min}}\tr(K-\tilde K)$$
    after using the cyclic invariance of the trace and the inequality
    $\tr(XY)\leq\lVert X\rVert\tr(Y)$
    for positive semidefinite matrices $X$ and $Y$.
    For matrix $B$, we define the diagonal matrix
    $\Delta=D^{-1/2}-\tilde D^{-1/2}$ and write
    $$B=D^{-1/2}\tilde K\Delta+\Delta\tilde KD^{-1/2}
        -\Delta\tilde K\Delta.$$
    Using the trace norm inequality
    $\lVert XYZ\rVert_*\leq\lVert X\rVert\lVert Y\rVert_*\lVert Z\rVert$
    together with $\lVert\tilde K\rVert_*=\tr(\tilde K)\leq \tr(K)$
    yields
    $$\lVert B\rVert_*\leq\bigl[\frac{2}{\sqrt{d_{\min}}}\lVert\Delta\rVert
        +\lVert\Delta\rVert^2\bigr]\tr(K).$$
    To bound $\norm{\Delta}=\max_i\abs{d_i^{-1/2}-\tilde d_i^{-1/2}}$
    we write
    $$\Bigl\lvert\frac{1}{\sqrt{d_i}}-\frac{1}{\sqrt{\tilde d_i}}\Bigr\rvert
        =\Bigl\lvert\frac{d_i-\tilde d_i}{\sqrt{d_i\tilde d_i}
        (\sqrt{d_i}+\sqrt{\tilde d_i})}\Bigr\rvert
        \leq\frac{\lvert d_i-\tilde d_i\rvert}{2\delta\sqrt{\delta}}$$
    using the inequalities $\sqrt{d_i\tilde d_i}\geq\delta$
    and $(\sqrt{d_i}+\sqrt{\tilde d_i})\geq 2\sqrt{\delta}$.
    Applying Lemma \ref{lem:dbound} yields
    $$\lVert\Delta\rVert\leq
        \frac{\sqrt{N}}{2\delta\sqrt{\delta}}\tr(K-\tilde K)$$
    and
    $$\lVert B\rVert_*\leq\Bigl[\frac{\sqrt{N}}{\delta^2}
        +\frac{N}{4\delta^3}\tr(K-\tilde K)\Bigr]\tr(K)\tr(K-\tilde K)$$
    after using $d_{\min}\geq\delta$.
    The triangle inequality
    $\lVert L-\tilde L\rVert_*\leq\lVert A\rVert_*+\lVert B\rVert_*$
    then yields the desired result.
\end{proof}

\begin{corollary2}[\ref{cor:symm-error2}]
    Under the assumptions of Theorem \ref{thm:rpc},
    choose $t\in\Rbb$ and the rank parameter $r$ such that
    \begin{equation}\label{eq:t-admit-copy}
        1<t<\frac{d_{\min}}{\sqrt{N}\zeta}.
    \end{equation}
    Then with probability at least $1-1/t$
    $$\lVert L-\tilde L\rVert_*<\Bigl[\frac{1}{d_{\min}}
        +\frac{\sqrt{N}\tr(K)}{(d_{\min}-\sqrt{N}t\zeta)^2}
        +\frac{N\tr(K)t\zeta}{4(d_{\min}-\sqrt{N}t\zeta)^3}\Bigr]t\zeta$$
    with $\tilde d_{\min}>d_{\min}-\sqrt{N}t\zeta>0$.
\end{corollary2}
\begin{proof}
    The random variable $\tr(K-\tilde K)$ is nonnegative, so the
    Markov inequality yields
    $$\Pbb\bigl\{\tr(K-\tilde K)\geq t\zeta\bigr\}\leq
        \frac{\Ebb\tr(K-\tilde K)}{t\zeta}\leq\frac{1}{t}$$
    using the fact that $\Ebb\tr(K-\tilde K)\leq\zeta$.
    For a value of $t\in\Rbb$ that satisfies the admissibility condition
    \eqref{eq:t-admit-copy},
    define the event $\Omega_t=\{\tr(K-\tilde K)<t\zeta\}$,
    so that $\Pbb(\Omega_t)\geq 1-1/t$.
    On the event $\Omega_t$, Lemma \ref{lem:dbound} yields
    $\lvert d_i-\tilde d_i\rvert<\sqrt{N}t\zeta$ for every $i$,
    which further implies
    $\tilde d_i>d_{\min}-\sqrt{N}t\zeta$ for every $i$.
    This leads to
    $$\tilde d_{\min}>d_{\min}-\sqrt{N}t\zeta>0.$$
    Based on the above,
    $\delta=\min(d_{\min},\tilde d_{\min})>d_{\min}-\sqrt{N}t\zeta$
    and $\tr(K-\tilde K)<t\zeta$, so Theorem \ref{thm:symm-error1}
    yields the desired result.
\end{proof}

\subsection{Bistochastic normalization}\label{app:proof-bistoch}
We prove the results of Section \ref{sec:error-bistoch}
related to the normalized matrices $P$ and $\tilde P$. 
We begin by proving the following lemma, which takes advantage
of the Gram structure of $P=G^\top G$
and $\tilde P=\tilde G^\top\tilde G$
introduced in \eqref{eq:gram} to essentially reduce a bound on
$P-\tilde P$ to one on $G-\tilde G$.
\begin{lemma}\label{lem:bistoch-gram}
    Under Assumptions \ref{asm:dpos} and \ref{asm:qpos},
    $$\norm{P-\tilde P}_*\leq(1+\norm{\tilde G})
        \norm{G-\tilde G}_*.$$
\end{lemma}
\begin{proof}
    We write
    $$P-\tilde P=G^\top G-\tilde G^\top\tilde G=
        G^\top(G-\tilde G)+(G-\tilde G)^\top\tilde G.$$
    Using the triangle inequality and the trace norm inequality
    $\norm{XY}_*\leq\norm{X}\norm{Y}_*$,
    \begin{equation*}
    \begin{split}
        \norm{P-\tilde P}_*&\leq\norm{G}\norm{G-\tilde G}_*
            +\norm{G-\tilde G}_*\norm{\tilde G}\\
        &\leq(1+\norm{\tilde G})\norm{G-\tilde G}_*.
    \end{split}
    \end{equation*}
For the last step, since $P$ is bistochastic, the Perron--Frobenius
theorem implies that $\norm{P}=\norm{G}^2=1$.
In contrast, even though the rows and columns of $\tilde P$
sum to 1, some of its entries may be negative, which means that
$\tilde P$ is generally not a bistochastic matrix.
As a result, the operator norm of $\tilde G$ may in general differ from 1.
\end{proof}

In addition to Lemma \ref{lem:dbound}, we are going to use
the following result.
\begin{lemma}\label{lem:qbound}
    Under Assumption \ref{asm:dpos}, we have
    $$\max_{0\leq i<N}\abs{q_i-\tilde q_i}\leq
        \frac{\sqrt{N}}{\delta_d}\bigl(1+\frac{\norm{K}}{\delta_d}\bigr)
        \tr(K-\tilde K).$$
\end{lemma}
\begin{proof}
    We write
    $$KD^{-1}-\tilde K\tilde D^{-1}=(K-\tilde K)D^{-1}
        +\tilde K(D^{-1}-\tilde D^{-1})$$
    and
    $$q-\tilde q=(K-\tilde K)D^{-1}1_N
        +\tilde K(D^{-1}-\tilde D^{-1})1_N$$
    where $q$ and $\tilde q$ denote the respective diagonals of
    $Q$ and $\tilde Q$.
    In what follows, we bound each term separately using the
    inequality $\norm{\cdot}_\infty\leq\norm{\cdot}_2$.
    We use the estimates $\norm{K-\tilde K}\leq\tr(K-\tilde K)$
    and $\norm{D^{-1}1_N}_2^2=\sum_{i=0}^{N-1}1/d_i^2\leq N/\delta_d^2$
    to bound the first term
    $$\norm{(K-\tilde K)D^{-1}1_N}_\infty\leq\norm{K-\tilde K}
        \norm{D^{-1}1_N}_2\leq\frac{\sqrt{N}}{\delta_d}\tr(K-\tilde K).$$
    For the second term, we write
    $\abs{d_i^{-1}-\tilde d_i^{-1}}=\abs{(\tilde d_i-d_i)/(d_i\tilde d_i)}
    \leq\abs{\tilde d_i-d_i}/\delta_d^2$
    and apply Lemma \ref{lem:dbound} to conclude that
    $$\sum_{i=0}^{N-1}\abs{d_i^{-1}-\tilde d_i^{-1}}^2
        \leq\frac{1}{\delta_d^4}
        \sum_{i=0}^{N-1}\abs{\tilde d_i-d_i}^2
        \leq\frac{N}{\delta_d^4}[\tr(K-\tilde K)]^2$$
    and
    $$\norm{\tilde K(D^{-1}-\tilde D^{-1})1_N}_\infty
        \leq\norm{\tilde K}\norm{(D^{-1}-\tilde D^{-1})1_N}_2
        \leq\norm{K}\frac{\sqrt{N}}{\delta_d^2}\tr(K-\tilde K).$$
\end{proof}

\begin{theorem2}[\ref{thm:bistoch-error1}]
    Under Assumptions \ref{asm:dpos} and \ref{asm:qpos},
    the normalized matrices $P$ and $\tilde P$ satisfy
    $$\norm{P-\tilde P}_*\leq(1+\norm{\tilde G})\Bigl[
        \frac{1}{\sqrt{q_{\min}}d_{\min}}
        +\frac{\sqrt{N}\tr(K)}{\sqrt{\delta_q}\delta_d^2}
        +\frac{\sqrt{N}\tr(K)}{2\delta_q^{3/2}\delta_d^2}
        \bigl(1+\frac{\norm{K}}{\delta_d}\bigr)\Bigr]
        \tr(K-\tilde K).$$
\end{theorem2}
\begin{proof}
    By Lemma \ref{lem:bistoch-gram}, it is sufficient to bound
    $G-\tilde G$ in the trace norm.
    To that end, we write
    \begin{equation*}
    \begin{split}
        G-\tilde G&=\tilde Q^{-1/2}\tilde K(D^{-1}-\tilde D^{-1})
            +Q^{-1/2}(K-\tilde K)D^{-1}\\
        &\quad+(Q^{-1/2}-\tilde Q^{-1/2})\tilde KD^{-1}
    \end{split}
    \end{equation*}
    and bound each term separately using the trace norm inequality
    $\norm{XYZ}_*\leq\norm{X}\norm{Y}_*\norm{Z}$.
    For the second term,
    $$\norm{Q^{-1/2}(K-\tilde K)D^{-1}}_*
        \leq\norm{Q^{-1/2}}\norm{K-\tilde K}_*\norm{D^{-1}}
        \leq\frac{\tr(K-\tilde K)}{\sqrt{q_{\min}}d_{\min}}.$$
    For the first term,
    $$\norm{\tilde Q^{-1/2}\tilde K(D^{-1}-\tilde D^{-1})}_*
        \leq\norm{\tilde Q^{-1/2}}\norm{\tilde K}_*
        \norm{D^{-1}-\tilde D^{-1}}
        \leq\frac{\sqrt{N}\tr(K)}{\sqrt{\delta_q}\delta_d^2}
        \tr(K-\tilde K)$$
    where we used Lemma \ref{lem:dbound} to obtain
    $$\norm{D^{-1}-\tilde D^{-1}}
        =\max_{0\leq i<N}\frac{\abs{d_i-\tilde d_i}}{d_i\tilde d_i}
        \leq\frac{1}{\delta_d^2}\max_{0\leq i<N}\abs{d_i-\tilde d_i}
        \leq\frac{\sqrt{N}}{\delta_d^2}\tr(K-\tilde K).$$
    Similarly, for the third term,
    \begin{equation*}
    \begin{split}
        \norm{(Q^{-1/2}-\tilde Q^{-1/2})\tilde KD^{-1}}_*
            &\leq\norm{(Q^{-1/2}-\tilde Q^{-1/2})}\norm{\tilde K}_*
            \norm{D^{-1}}\\
            &\leq\frac{\sqrt{N}\tr(K)}{2\delta_q^{3/2}\delta_d^2}
            \bigl(1+\frac{\norm{K}}{\delta_d}\bigr)\tr(K-\tilde K)
    \end{split}
    \end{equation*}
    where we used Lemma \ref{lem:qbound} to obtain
    \begin{equation*}
    \begin{split}
        \norm{(Q^{-1/2}-\tilde Q^{-1/2})}
            =\max_{0\leq i<N}\abs{\frac{1}{\sqrt{q_i}}
            -\frac{1}{\sqrt{\tilde q_i}}}
            &=\max_{0\leq i<N}\abs{\frac{\tilde q_i-q_i}
            {\sqrt{q_i\tilde q_i}(\sqrt{q_i}+\sqrt{\tilde q_i})}}\\
            &\leq\frac{1}{2\delta_q^{3/2}}
            \max_{0\leq i<N}\abs{\tilde q_i-q_i}\\
            &\leq\frac{\sqrt{N}}{2\delta_q^{3/2}\delta_d}
            \bigl(1+\frac{\norm{K}}{\delta_d}\bigr)\tr(K-\tilde K).
    \end{split}
    \end{equation*}
    An application of Lemma \ref{lem:bistoch-gram}
    and the triangle inequality then yields the desired result.
\end{proof}

\begin{corollary2}[\ref{cor:bistoch-error2}]
    Under the assumptions of Theorem \ref{thm:rpc},
    let $t\in\Rbb$ and
    $$\alpha_t=d_{\min}-\sqrt{N}t\zeta\qquad
        \beta_t=\frac{\sqrt{N}t\zeta}{\alpha_t}
        \bigl(1+\frac{\norm{K}}{\alpha_t}\bigr)$$
    and choose $t$ and the rank parameter $r$ such that
    \begin{equation}\label{eq:t-admit2-copy}
        t>1\qquad\alpha_t>0\qquad\beta_t<q_{\min}.
    \end{equation}
    Then with probability at least $1-1/t$
    \begin{equation*}
    \begin{split}
        \norm{P-\tilde P}_*<(1+\norm{\tilde G})\Bigl[
            &\frac{1}{\sqrt{q_{\min}-\beta_t}\alpha_t}
            +\frac{\sqrt{N}\tr(K)}{\sqrt{q_{\min}-\beta_t}\alpha_t^2}\\
            &+\frac{\sqrt{N}\tr(K)}{2(q_{\min}-\beta_t)^{3/2}\alpha_t^2}
            \bigl(1+\frac{\norm{K}}{\alpha_t}\bigr)\Bigr]t\zeta
    \end{split}
    \end{equation*}
    with $\tilde d_{\min}>\alpha_t>0$
    and $\tilde q_{\min}>q_{\min}-\beta_t>0$.
\end{corollary2}
\begin{proof}
    We proceed again by constructing a high probability event
    $\Omega_t$ and conditioning on it to produce a deterministic bound,
    mirroring the approach used in proving
    Corollary \ref{cor:symm-error2}.
    Since the random variable $\tr(K-\tilde K)$ is nonnegative,
    the Markov inequality yields
    $$\Pbb\{\tr(K-\tilde K)\geq t\zeta\}
        \leq\frac{\Ebb\tr(K-\tilde K)}{t\zeta}\leq\frac{1}{t}.$$
    For $t\in\Rbb$ satisfying the admissibility condition
    \eqref{eq:t-admit2-copy}, define the event
    $\Omega_t=\{\tr(K-\tilde K)<t\zeta\}$, so that
    $\Pbb(\Omega_t)\geq 1-1/t$.
    On the event $\Omega_t$, Lemma \ref{lem:dbound}
    yields $\abs{d_i-\tilde d_i}<\sqrt{N}t\zeta$ for every $i$;
    this implies $\tilde d_i>d_{\min}-\sqrt{N}t\zeta$ for every $i$,
    and $\tilde d_{\min}>\alpha_t>0$.
    As a result, $\delta_d>\alpha_t>0$ on $\Omega_t$ for $t$
    satisfying \eqref{eq:t-admit2-copy}.
    On the same event $\Omega_t$, Lemma \ref{lem:qbound}
    can now be used to produce 
    $$\max_i\abs{q_i-\tilde q_i}<\frac{\sqrt{N}}{\delta_d}
        \bigl(1+\frac{\norm{K}}{\delta_d}\bigr)t\zeta<\beta_t.$$
    A similar argument as above shows that
    $\tilde q_i>q_{\min}-\beta_t$ for every $i$,
    which leads to $\tilde q_{\min}>q_{\min}-\beta_t>0$.
    Namely, $\delta_q>q_{\min}-\beta_t>0$ on $\Omega_t$
    for $t$ satisfying \eqref{eq:t-admit2-copy}.
    An application of Theorem \ref{thm:bistoch-error1}
    then yields the desired result.
\end{proof}

\section*{Use of Generative AI}
We used the Claude Opus family of AI models to assist us in the development
of the error bounds and of the code used for our numerical results.
The full manuscript text, including the error bounds and their proofs,
was written by the authors.

\section*{Acknowledgments}
DG acknowledges support from the US Department of Energy under
grant DE-SC0025101.
CV was supported as a postdoctoral researcher from this grant.
We thank the reviewers for their careful and constructive feedback,
which prompted us to expand the scope of this work and improve its
presentation.

\medskip

Received xxxx 20xx; revised xxxx 20xx; early access xxxx 20xx.
\medskip
\end{document}